\documentclass[11pt]{article}
\usepackage{amssymb}
\usepackage{amsmath}
\usepackage{mathrsfs}
\usepackage{graphics}
\usepackage{graphicx}
\usepackage{subfigure}
\usepackage[T1]{fontenc}
\usepackage{latexsym,amssymb,amsmath,amsfonts,amsthm}\usepackage{txfonts}
\newcommand{\R}{{\mathbb R}}

\newcommand{\Sp}{{\mathbb S}}

\newcommand{\no}{\nonumber}
\newcommand{\be}{\begin{eqnarray}}
\newcommand{\ben}{\begin{eqnarray*}}
\newcommand{\en}{\end{eqnarray}}
\newcommand{\enn}{\end{eqnarray*}}
\newcommand{\ba}{\backslash}
\newcommand{\pa}{\partial}

\newcommand{\ov}{\overline}
\newcommand{\I}{{\rm Im}}

\newcommand{\g}{\gamma}

\newcommand{\Om}{\Omega}

\newcommand{\al}{\alpha}
\newcommand{\la}{\lambda}

\newcommand{\ol}{\overline}

\newcommand{\half}{\frac{1}{2}}
\newcommand{\md}{\mathrm{d}}
\newcommand{\na}{\nabla}
\newtheorem{theorem}{Theorem}[section]

\newtheorem{remark}[theorem]{Remark}

\newtheorem{algorithm}{Algorithm}[section]

\begin{document}
\renewcommand{\theequation}{\arabic{section}.\arabic{equation}}
\title{\bf A Newton method for simultaneous reconstruction of an interface
and a buried obstacle from far-field data}
\author{Haiwen Zhang,\ \ Bo Zhang\\
LSEC and Institute of Applied Mathematics\\
Academy of Mathematics and Systems Science\\
Chinese Academy of Sciences\\
Beijing 100190, China\\
({\sf zhanghaiwen@amss.ac.cn},\ {\sf b.zhang@amt.ac.cn})
}
\date{}
\maketitle

\begin{abstract}
This paper is concerned with the inverse problem of scattering of time-harmonic
acoustic waves from a penetrable and buried obstacles.
By introducing a related transmission scattering problem, a Newton iteration method
is proposed to simultaneously reconstruct both the penetrable interface and the buried
obstacle inside from far-field data. 
A main feature of our method is that we do not need to know the type of boundary conditions
on the buried obstacle. In particular, the boundary condition on the buried obstacle can also be 
determined simultaneously by the method. Finally, numerical examples using multi-frequency data 
are carried out to illustrate the effectiveness of our method.
\end{abstract}

\section{Introduction}\label{sec1}

In this paper, we consider the problem of scattering of time-harmonic acoustic
waves from a penetrable obstacle with a buried impenetrable obstacle inside.
Such problems occur in many applications such as radar, remote sensing, geophysics and
nondestructive testing.

Let $S_0\in C^2$ denote a simple closed smooth curve in $\R^2$.
Then $\mathbb{R}^2$ is divided into the unbounded part $\Om_0$ and the bounded
part $\Om$ by $S_0$. We assume that $\Om_0$ is filled with a homogeneous medium with
the constant refractive index $1$ and $\Om$ is filled with another
homogeneous medium with the constant refractive index $n_1>0$.
Assume further that $\Om_2\subset\subset\Om$ is an impenetrable obstacle with a $C^2$
boundary $S_1$ and let $\Om_1=\Om\ba\ol{\Om}_2$.

\begin{figure}[htbp]
\centering
\includegraphics[width=4in]{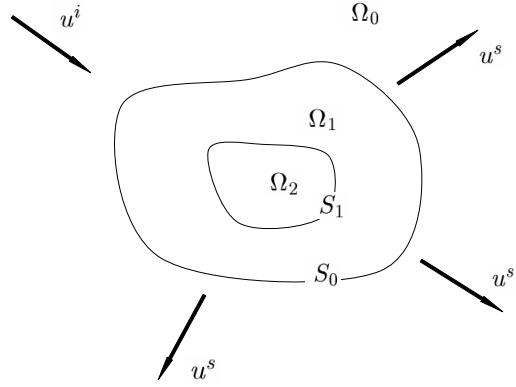}
\caption{Scattering by a penetrable and a buried obstacles}
\label{pic1}  
\end{figure}

Suppose an incident plane wave $u^i(x,d)=e^{ik_0x\cdot d}$,
where $k_0=\omega/c>0$ is the wave number with $\omega$ and $c$ being the wave
frequency and speed in $\Omega_0$ and $d\in S^1$ is the incident direction,
is incident on the penetrable obstacle $\Om$ from the unbounded domain $\Om_0$.
Then the total field $u=u^i+u^s$, which is the sum of the incident wave $u^i$ and
the scattered wave $u^s$ in $\Om_0$, satisfies the following Helmholtz equations
together with the transmission condition on the interface $S_0$ and
the boundary condition on the boundary $S_1$:
\be\label{eq16}
\Delta u+k^2_0 u=0&&in\;\;\Omega_0\\ \label{eq17}
\Delta u+k^2_1 u=0&&in\;\;\Omega_1\\ \label{eq14}
u_+=u_-,\quad\frac{\pa u_+}{\pa\nu}=\la_0\frac{\pa u_-}{\pa\nu}&&on\;\;S_0\\ \label{eq15}
\mathscr{B}(u)=0&&on\;\;S_1\\ \label{rc}
\lim_{r\to\infty}\sqrt{r}\left(\frac{\pa u^s}{\pa r}-ik_0u^s\right)=0,&&\quad r=|x|
\en
where $k_1>0$ is the wave number in $\Om_1$ given by $k_1^2=k^2_0 n_1$,
$\la_0$ is a positive constant depending on the property of the medium in $\Om$
and $\Om_0$, $\nu$ is the exterior normal vector on $S_0$ and +/- denote
the limit from the exterior and interior of the boundary, respectively.
The boundary condition $\mathscr{B}(u)=0$ represents the Dirichlet or
Neumann boundary condition depending on the physical property of the obstacle $\Om_2$.
It is well-known that the Sommerfield radiation condition (\ref{rc}) implies
the following asymptotic behavior of the scattered field $u^s$ (see \cite{ColtonKress98}):
\ben
u^s(x,d,k_0)=\frac{e^{ik_0r}}{\sqrt{r}}\left(u_\infty(\hat{x},d,k_0)
+O(\frac{1}{r})\right),\qquad r=|x|\to\infty,
\enn
where $\hat{x}=x/r$ is the observation direction and $u_\infty$ is called the far
field patten of the scattered field $u^s$.

By using a variational or integral equation method it has been proved
in \cite{LiuZhang2010,LiuZhangHu2010} that the scattering problem
(\ref{eq16})-(\ref{rc}) has a unique solution.
It was further proved in \cite{LiuZhang2010} that
$S_0,S_1$ and the boundary condition on $S_1$ can be uniquely
determined by the far-field pattern $u_\infty(\hat{x},d,k_0)$ for all $\hat{x},d\in S^1$
and a fixed $k_0$ if $\lambda_0\neq 1$ is known.
In this paper, we are interested in the numerical reconstruction of the interface $S_0$
and the buried impenetrable obstacle $S_1$ together with the boundary condition on $S_1$
from the far-field pattern $u_\infty(\hat{x},d,k_0)$, given $\la_0>0$ and
the refractive index $n_1>0$.

Many reconstruction algorithms have been developed for the numerical reconstruction of
the interface $S_0$ in the case when the buried obstacle $\Om_0=\emptyset$
(see, e.g., the iteration methods in \cite{AK,HS} and the singular sources method
in \cite{PS}, the linear sampling method in \cite{M05} and the factorization method
in \cite{Kirsch1,Kirsch2} for the case when $\la_0=1$).
In the case when the interface $S_0$ is known and $\la_0=1$, a reciprocity gap functional
method was introduced in \cite{CFH,CH05,DS06,DS07} to reconstruct the buried obstacle $\Om_0$
from near field Cauchy data, whilst in \cite{Yaman} an iteration method was proposed
to recover the buried sound-soft obstacle from the far field data.
Recently, for the case when $\Om_0\not=\emptyset$, Cakoni et al. \cite{CDS} introduced a
multistep reciprocity gap functional method for reconstructing the interface $S_0$ and
the buried obstacle $\Om_0$ recursively from near field Cauchy data under the condition
that the boundary condition on $S_1$ or $\Om_0$ is known in advance,
whilst in \cite{YZZ} the factorization method was extended to reconstruct
the interface $S_0$ from the far field data, without knowing the buried obstacle $\Om_0$.

In this paper we develop a Newton-type iterative method to simultaneously reconstruct the
interface $S_0$ and the impenetrable obstacle $S_1$ together with the boundary
condition on $S_1$ from the far-field pattern $u_\infty(\cdot,d,k_0)$ for a finite number
of incident directions $d=d_j,$ $j=1,...,P$ and a finite number of frequencies or wavenumbers
$k_0=k_{0l},\;l=1,...,Q,$ given $\la_0>0$ and the refractive index $n_1>0$.
A main feature of our method is that the boundary condition on the buried obstacle $S_1$
does not need to know in advance; in fact, the boundary condition on $S_1$ can be
determined by the iterative method. This is different from the previous iteration methods
which require to know the boundary condition on the obstacle to be reconstructed
as a priori information.
To simultaneously reconstruct $S_0, S_1$ and the boundary condition on $S_1$,
we first introduce a related transmission scattering problem by replacing the impenetrable
obstacle $\Om_0$ with a penetrable one with an arbitrarily chosen wavenumber $k_2$
and having a transmission condition on $S_1$ with the unknown transmission constant $\la_1$.
We then use a Newton method to simultaneously reconstruct $S_0, S_1$ and $\la_1$.
The boundary condition on $S_1$ can thus be determined in terms of the value of the
reconstructed $\la_1$. See Section \ref{sec2} for detailed discussion.

Over the last 20 years, many work has been done related to Newton-type methods.
For example, \cite{Kirsch93} gives the first rigorous characterization of Frechet derivatives 
of the far-field operator corresponding to an impenetrable, sound-soft obstacle embedded in a 
homogeneous medium by employing a variational method, and \cite{Potthast1994} obtained the same 
result via the integral equation method. For the other cases such as impedance and transmission 
problems, we refer to \cite{Hettlich1995,HS}. 
It should be pointed out that the convergence of Newton-type methods for inverse scattering
is far from complete. For partial results, the reader is referred to 
\cite{HS,Hohage1998,Hohage99,Potthast01}. 
In addition, much effort has been made to reduce the computation cost of the Newton-type methods,
such as the hybrid method \cite{KressSerranho2005,KressSerranho2007} and 
the methods of nonlinear integral equations \cite{AK,KressRundell2005}.

This paper is organized as follows. In Section \ref{sec2}, a related transmission 
scattering problem is presented. In Section \ref{sec3} the integral equation method
is proposed for solving the related transmission problem, and some quadrature rules are
discussed briefly. The Newton method is presented in Section \ref{sec4},
and numerical examples using multi-frequency data are given in Section \ref{sec5} to illustrate 
the effectiveness of our method. In the appendix, we give characterizations of the Frechet 
derivatives of the far field operators for the related transmission scattering problem in 
Section \ref{sec2}, employing the variational method.

\section{A related transmission scattering problem}\label{sec2}

As discussed in the previous section, to simultaneously reconstruct $S_0$,
$S_1$ and the boundary condition on $S_1$, we consider the following transmission
scattering problem which is related to our original scattering problem:
find the total wave $u\coloneqq u^i+u^s$ (which is the sum of the incident field
$u^i$ and the scattered field $u^s$ in $\Om_0$) such that
\be\label{eq20}
\Delta u+k^2_0 u=0& \mbox{in}\;\;\Omega_0\\ \label{eq21}
\Delta u+k^2_1 u=0& \mbox{in}\;\;\Omega_1\\ \label{eq22}
\Delta u+k^2_2 u=0& \mbox{in}\;\;\Omega_2\\ \label{eq23}
u_+=u_-,\quad\frac{\pa u_+}{\pa\nu}=\la_0\frac{\pa u_-}{\pa\nu}&\mbox{on}\;\;S_0 \\ \label{eq24}
u_+=u_-,\quad\frac{\pa u_+}{\pa\nu}=\la_1\frac{\pa u_-}{\pa\nu}&\mbox{on}\;\;S_1
\en
where $u^i,k_0,k_1,\la_0$ are as defined in Section \ref{sec1},
$k_2$ is an arbitrarily chosen positive constant, $\la_1$ is a constant to be determined,
and $u^s$ satisfies the Sommerfeld radiation condition (\ref{rc}).
This problem is obtained from the scattering problem (\ref{eq16})-(\ref{rc}) by regarding
the impenetrable obstacle $\Om_0$ as a penetrable one having an artificial wavenumber
$k_2$ and a transmission condition on its boundary $S_1$ with the unknown
transmission constant $\la_1$.

We observe that the Dirichlet and Neumann boundary conditions can be considered as the limiting
cases $\la_1\to\infty$ and $\la_1\to0$, respectively, of the transmission condition (\ref{eq24}).
Based on this observation, we will develop a Newton method to simultaneously reconstruct $S_0,S_1$
and $\la_1$ and then use the value of the reconstructed $\la_1$ to determine the type of
boundary conditions on the boundary $S_1$ of the impenetrable obstacle $\Om_0$.
Precisely, if the value of the reconstructed $\la_1$ is very close to zero, then the boundary
condition on $S_1$ will be regarded as a Neumann condition. If the value of the reconstructed $\la_1$
is far away from zero, we define $\tau_1\coloneqq 1/\la_1$, $U\coloneqq\la_1 u|_{\Om_2}$
and replace (\ref{eq22}) and (\ref{eq24}), respectively, by
\be\label{eq22a}
\Delta U+k^2_2 U=0&\mbox{in}\;\;\Omega_2\\ \label{eq24a}
u_+=\tau_1 U_-,\quad\frac{\pa u_+}{\pa\nu}=\frac{\pa U_-}{\pa\nu}&\mbox{on}\;\;S_1.
\en
We then use the Newton method to simultaneously reconstruct $S_0, S_1$ and $\tau_1$.
If the value of the reconstructed $\tau_1$ is very close to zero, then the boundary condition
on $S_1$ will be regarded as a Dirichlet condition.

\begin{remark}\label{re1} {\rm
During the preparation of this paper we were notified by the authors of \cite{LiuShang}
that they have recently also studied the asymptotic behavior of the solution to the scattering
problem by a lossy inhomogeneous medium buried in a layered homogeneous medium as the medium
mass density is very large in the context of cloaking or invisibility; precisely, they
proved mathematically that the solution of the transmission problem (\ref{eq20})-(\ref{eq24})
and (\ref{rc}) converges to the solution of the scattering problem (\ref{eq16})-(\ref{rc})
with a Neumann boundary condition $\mathscr{B}(u):=\pa u/\pa\nu=0$ on $S_1$ as $\la_1\to0+$
in the case when $\I(k^2_2)>0$.
}
\end{remark}

Based on the above observation, we assume that $k_0$, $k_1$, $k_2$, $\la_0$ are known and
introduce the far field operator $F$ which maps a set of admissible boundaries $S_0$ and $S_1$
and the constant $\la_1$ into the far field pattern of the scattered solution
to the corresponding transmission scattering problem (\ref{eq20})-(\ref{eq24}) and (\ref{rc}).
We also introduce the far field operator $\tilde{F}$ which maps a set of admissible
boundaries $S_0$ and $S_1$ and the constant $\tau_1$ into the far field pattern of the solution
to the transmission scattering problem (\ref{eq20})-(\ref{eq21}), (\ref{eq23}),
(\ref{eq22a})-(\ref{eq24a}) and (\ref{rc}).
We will solve the equation
\be\label{eq19}
F(S_0,S_1,\la_1)=u_\infty
\en
for $S_0,S_1,\la_1$; as discussed above, for the case when $\la_1$ is very large, we
solve the equation
\be\label{eq25}
\tilde{F}(S_0,S_1,\tau_1)=u_\infty
\en
for $S_0,S_1,\tau_1$. Here, the value of $\la_1$ or $\tau_1$ reveal the type of
boundary conditions on $S_1$ and $u_\infty$ is the far field pattern of the original
scattering problem (\ref{eq16})-(\ref{rc}).

Following the idea of \cite{Kirsch93,Hettlich1995}, we can prove that if $S_0$ and $S_1$
are $C^2$ boundary then $F$ (or $\tilde{F}$) is Frechet differentiable
with respect to $S_0$, $S_1$ and $\la_1$ (or $\tau_1$) for the case when
$\la_1\in(0,\infty)$ (or $\tau\in(0,\infty)$) (see Appendix for details).
Motivated by this we propose a Newton iteration method to recover the boundaries $S_0$, $S_1$
and the boundary condition on $S_1$ in Section \ref{sec4}.


\section{The integral equation method for the direct problem}\label{sec3}
\setcounter{equation}{0}

For the Newton method, we need to compute the Frechet derivative of $F$ or $\tilde{F}$
which is characterized by the solution
of the problem (\ref{eq20})-(\ref{eq24}) and (\ref{rc})
with certain inhomogeneous boundary data (see Appendix).
For convenience we consider the following more general transmission problem:
find $u^s\in C^2(\Om_0)\bigcap C^{1,\sigma}(\ov{\Om}_0)$ and
$u\in C^2(\Om_1)\bigcap C^{1,\sigma}(\ov{\Om}_1)\bigcap C^2(\Om_2)
\bigcap C^{1,\sigma}(\ov{\Om}_2)$ with $0<\sigma<1$ such that
\be\label{eq39}
\Delta u^s+k^2_0 u^s=0&\mbox{in}\;\;\Omega_0\\ \label{eq40}
\Delta u+k^2_1 u=0&\mbox{in}\;\;\Omega_1\\ \label{eq10}
\Delta u+k^2_2 u=0&\mbox{in}\;\;\Omega_2\\ \label{eq41}
u^s_+-u_-=f_1,\quad
\frac{\pa u^s_+}{\pa\nu}-\la_0\frac{\pa u_-}{\pa\nu}=f_2&\mbox{on}\;\;S_0\\ \label{eq11}
u_+-u_-=f_3,\quad\frac{\pa u_+}{\pa\nu}-\la_1\frac{\pa u_-}{\pa\nu}=f_4&\mbox{on}\;\;S_1
\en
and $u^s$ satisfies the radiation condition (\ref{rc}),
where $(f_1,f_2,f_3,f_4)\in X\coloneqq C^{1,\sigma}(S_0)\times C^{0,\sigma}(S_0)
\times C^{1,\sigma}(S_1)\times C^{0,\sigma}(S_1)$.

\begin{remark}\label{re2} {\rm
For our scattering problem (\ref{eq20})-(\ref{eq24}) and
(\ref{rc}), $(f_1,f_2,f_3,f_4)=(-u^i,-{\pa u^i}/{\pa\nu},0,0)$.
}
\end{remark}

The transmission problem (\ref{eq39})-(\ref{eq11}) and (\ref{rc})
has been studied in \cite{LZ12}, employing the integral equation method.
In this paper we will also use the integral equation method to deal with the inverse problem.
To this end, for $l=0,1,2$ let $\Phi_l(x,y)$ denote the fundamental solution of the Helmholtz
equation $\Delta +k^2_l$ in $\R^2$ and, for $i,j=0,1$ define
the boundary integral operators $S_{ijl}$, $K_{ijl}$, $K^T_{ijl}$, $T_{ijl}$ by
\ben
(S_{ijl}\psi)(x)=\int_{S_j}\Phi_l(x,y)\psi(y)ds(y)&&x\in S_i\\
(K_{ijl}\psi)(x)=\int_{S_j}\frac{\pa\Phi_l(x,y)}{\pa\nu(y)}\psi(y)ds(y)&&x\in S_i\\
(K^T_{ijl}\psi)(x)=\int_{S_j}\frac{\pa\Phi_l(x,y)}{\pa\nu(x)}\psi(y)ds(y)&&x\in S_i\\
(T_{ijl}\psi)(x)=\frac{\pa}{\pa\nu(x)}\int_{S_j}\frac{\pa\Phi_l(x,y)}{\pa\nu(y)}\psi(y)ds(y)&&x\in S_i
\enn
We also define the far field operators $S_\infty$ and $K_\infty$
corresponding to $S_{ij0}$, $K_{ij0}$ by
\ben
(S_\infty\psi)(\hat{x})=\frac{e^{i\pi/4}}{\sqrt{8\pi k_0}}
  \int_{S_0}e^{-ik_0\hat{x}\cdot y}\psi(y)ds(y) &&\hat{x}\in\Sp^1\\
(K_\infty\psi)(\hat{x})=\frac{e^{i\pi/4}}{\sqrt{8\pi k_0}}
  \int_{S_0}\frac{\pa e^{-ik_0\hat{x}\cdot y}}{\pa\nu(y)}\psi(y)ds(y) &&\hat{x}\in\Sp^1
\enn
where $\Sp^1$ is the unit sphere in $\R^2$.
The reader is referred to \cite{ColtonKress98} for mapping properties of these operators
in the spaces of continuous and H\"{o}lder continuous functions.

Following \cite{LZ12}, we seek a solution to the transmission problem (\ref{eq39})-(\ref{eq11})
and (\ref{rc}) in the following form:
\ben
u^s(x)&=&\int_{S_0}\left[\la_0\frac{\pa\Phi_0(x,y)}{\pa\nu(y)}\psi_1(y)+\Phi_0(x,y)\psi_2(y)\right]ds(y)
   \qquad x\in\Om_0\\
u(x)&=&\int_{S_0}\left[\frac{\pa\Phi_1(x,y)}{\pa\nu(y)}\psi_1+\Phi_1(x,y)\psi_2(y)\right]ds(x)\\
  &&+\int_{S_1}\left[\la_1\frac{\pa\Phi_1(x,y)}{\pa\nu(y)}\psi_3(y)+\Phi_1(x,y)\psi_4(y)\right]ds(y)
  \qquad x\in\Om_1\\
u(x)&=&\int_{S_1}\left[\frac{\pa\Phi_2(x,y)}{\pa\nu(y)}\psi_3+\Phi_2(x,y)\psi_4(y)\right]ds(x)
\qquad x\in\Om_2
\enn
with $(\psi_1,\psi_2,\psi_3,\psi_4)\in X$.
From the jump relations of the potential operators, it is known that the potentials $u^s$ and $u$
defined above solve the problem (\ref{eq39})-(\ref{eq11}) and (\ref{rc}) if
$\Psi=(\psi_1,\psi_2,\psi_3,\psi_4)^T$ solves the operator equation
\be\label{eq33}
(I+A)\Psi=R
\en
where $R=(\mu_0f_1,-\mu_0 f_2,\mu_1f_3,-\mu_1f_4)^T$, $\mu_i=2/(1+\la_i)$, i=0, 1, $I$ is an
identity operator and $A:X\rightarrow X$ is given by
\ben
\left(
\begin{array}{cccc}
\mu_0(\la_0 K_{000}-K_{001}) & \mu_0(S_{000}-S_{001}) & -\mu_0\la_1 K_{011} & -\mu_0 S_{011} \\
-\mu_0\la_0(T_{000}-T_{001})&-\mu_0(K^T_{000}-\la_0K^T_{001})&\mu_0\la_0\la_1T_{011}&\mu_0\la_0K^T_{011}\\
\mu_1 K_{101}& \mu_1 S_{101} & \mu_1(\la_1 K_{111}-K_{112}) & \mu_1(S_{111}-S_{112}) \\
-\mu_1 T_{101} & -\mu_1 K^T_{101} & -\mu_1\la_1(T_{111}-T_{112}) & -\mu_1(K^T_{111}-\la_1 K^T_{112} )\\
\end{array}
\right)
\enn
It was proved in \cite{LZ12} that the system (\ref{eq33}) of integral equations is uniquely
solvable in $X$ (see Theorem 2.4 in \cite{LZ12} for details).
Further, by the asymptotic behavior of the Green function \cite{ColtonKress98}, we have
\be\label{eq34}
u_\infty(\hat{x})=\la_0(S_\infty\psi_1)(\hat{x})+(K_\infty\psi_2)(\hat{x}),\qquad\hat{x}\in S^1,
\en
where $(\psi_1,\psi_2)\in C^{1,\sigma}(S_0)\times C^{0,\sigma}(S_0)$ are the first two components
of the unique solution to (\ref{eq33}).

To generate the far field data, we use the Nystr\"{o}m method to solve the system (\ref{eq33}) of
integral equations.
We assume that the boundaries $S_j,j=0,1$, can be parameterized
as $x(t)=(x_1(t),x_2(t))$ for $0\le t\le2\pi$, where $x(t)$ is an analytic and $2\pi$-periodic
function. Then, by the change of variables the integral operators $S_{ijl},K_{ijk},K^T_{ijl}$,
$i,j=0,1,l=0,1,2$, and $T_{ijl}$ ($i,j=0,1,l=0,1,2,i\neq j$) can be reduced to the form
\ben
\int^{2\pi}_0\left[k_1(t,\tau)\ln\left(4sin^2\frac{t-\tau}{2}\right)
+k_2(t,\tau)\right]\phi(\tau)d\tau
\enn
and $T_{ijl}$ ($i=j$) can be reduced to the form
\ben
\frac{1}{2\pi}\int^{2\pi}_0\cot\frac{\tau-t}{2}\phi'(\tau)\md\tau
+\int^{2\pi}_0\left[k_1(t,\tau)\ln\left(4sin^2\frac{t-\tau}{2}\right)
+k_2(t,\tau)\right]\phi(\tau)\md\tau
\enn
where $k_1(t,\tau),k_2(t,\tau)$ are analytic and $2\pi$-periodic functions.

In this paper we use the following quadrature rule: for $0\le t\le2\pi$,
\ben
\int^{2\pi}_{0}ln\left(4\sin^2\frac{t-\tau}{2}\right)f(\tau)d\tau
  &\approx&\sum^{2n-1}_{j=0}R^{(n)}_j(t)f(t_j)\\
\int^{2\pi}_{0}f(\tau)d\tau&\approx&\frac{\pi}{n}\sum^{2n-1}_{j=0}f(t_j)\\
\frac{1}{2\pi}\int^{2\pi}_0\cot\frac{\tau-t}{2}f'(\tau)d\tau
  &\approx&\sum^{2n-1}_{j=0}T^{(n)}_j(t)f(t^{(n)}_j)
\enn
where, for $j=0,\ldots,2n-1,$
\ben
R^{(n)}_j(t)&\coloneqq& -\frac{2\pi}{n}\sum^{n-1}_{m=1}
  \frac{1}{m}\cos[m(t-t_j)]-\frac{\pi}{n^2}\cos[n(t-t^{(n)}_j)]\\
T^{(n)}_j(t)&\coloneqq& -\frac{1}{n}\sum^{n-1}_{m=1}m\cos[m(t-t^{(n)}_j)]
  -\half\cos[n(t-t^{(n)}_j)]
\enn
After finding the solution $\Psi=(\psi_1,\psi_2,\psi_3,\psi_4)^T$ to the system of integral
equations, (\ref{eq33}), we can use (\ref{eq34}) to compute the far field pattern
(cf. \cite{ColtonKress98,Kress95,Yaman}).

\section{The inverse problem}\label{sec4}
\setcounter{equation}{0}

In this section, we propose a Newton iteration method for the inverse problem.
We assume that $S_0$ and $S_1$ are starlike surfaces which can be parameterized
by $\g_0(\theta)$ and $\g_1(\theta)$, respectively:
\ben
\g_{0}(\theta)&=&(0,0)^T+r_0(\theta)\left(\cos(\theta),\sin(\theta)\right)^T\\
\g_{1}(\theta)&=&(a_1,a_2)^T+r_1(\theta)\left(\cos(\theta),\sin(\theta)\right)^T
\enn
where $\theta\in[0,2\pi]$ and $(0,0)$ and $(a_1,a_2)$ are their centers, respectively.
We will consider the noised perturbation $u^\delta_\infty$ of the far field data $u_\infty$
in the sense that $||u^\delta_\infty-u_\infty||_{L^2(S^1)}\le\delta||u_\infty||_{L^2(S^1)}$,
where $\delta\ge0$ is called the noise ratio.
Thus, given $u^\delta_\infty(\cdot,d_i,k_0),i=1,...,P$ for a fixed $k_0$, we
rewrite (\ref{eq19}) and (\ref{eq25}), respectively, as
\be\label{eq35}
F_i(\g_{0},\g_{1},\la_1)&\thickapprox& u^\delta_{\infty,i},\quad i=1,...,P\\ \label{eq36}
\tilde{F}_i(\g_{0},\g_{1},\tau_1)&\thickapprox& u^\delta_{\infty,i},\quad i=1,...,P
\en
where $u^\delta_{\infty,i}=u^\delta_{\infty}(\cdot,d_i,k_0)$ and $F_i,\tilde{F}_i$
represent the far field operators (\ref{eq35}), (\ref{eq36}) corresponding to the incident
plane wave with direction $d_i$.

We now derive the Newton method for (\ref{eq35}). The method is similar for (\ref{eq36}).
Assume that, after the $m$-th iteration, we have the approximate data $\g^{(m)}_0,\g^{(m)}_1,
\la^{(m)}_1$ with $\la^{(m)}_1\neq 0$, where
\ben
\g^{(m)}_{0}(\theta)&=&(0,0)^T+r^{(m)}_0(\theta)(\cos(\theta),\sin(\theta))^T\\
\g^{(m)}_{1}(\theta)&=&(a^{(m)}_1,a^{(m)}_2)^T+r^{(m)}_1(\theta)(\cos(\theta),\sin(\theta))^T
\enn
Then the linearized form of (\ref{eq35}) becomes
\be\label{eq-l}
\sum^{1}_{l=0}\left.\frac{\pa F_i}{\pa\g_l}\right|_{(\g^{(m)}_0,\g^{(m)}_1,\la^{(m)}_1)}
(\triangle\g^{(m)}_l)+\left.\frac{\pa F_i}{\pa\la_1}
\right|_{(\g^{(m)}_0,\g^{(m)}_1,\la^{(m)}_1)}(\triangle\la^{(m)}_1)
+F_i(\g^{(m)}_0,\g^{(m)}_1,\la^{(m)}_1)\thickapprox u^\delta_{\infty,i}
\en
for $i=1,...,P,$ where $\Delta\g^{(m)}_l(\theta)=\g^{(m+1)}_l(\theta)-\g^{(m)}_l(\theta)$, $l=0,1,$
and $\Delta\la^{(m)}_1=\la^{(m+1)}_1-\la^{(m)}_1$ are the updates to be determined.
The linearized equation (\ref{eq-l}) is ill-posed, so a regularization strategy is needed.
Following \cite{Hohage99}, we use the Levenberg-Marquardt algorithm, that is, we seek a solution
$(\Delta\g^{(m)}_0,\Delta\g^{(m)}_1,\Delta\la^{(m)}_1)$ with $\Delta\la^{(m)}_1\in\R$ and
\ben
\Delta\gamma^{(m)}_{0}(\theta)&=&(0,0)^T+\Delta r^{(m)}_0(\theta)(\cos(\theta),\sin(\theta))^T\\
\Delta\gamma^{(m)}_{1}(\theta)&=&(\Delta a^{(m)}_1,\Delta a^{(m)}_2)^T
+\Delta r^{(m)}_1(\theta)(\cos(\theta),\sin(\theta))^T,
\enn
where $\Delta r^{(m)}_0,\Delta r^{(m)}_1\in H^s(0,2\pi),$ $s\geq0,$ and
$\Delta a^{(m)}_1,\Delta a^{(m)}_2\in\R$, such that
$(\Delta\g^{(m)}_0,\Delta\g^{(m)}_1,\Delta\la^{(m)}_1)$ is a solution to the minimization problem:
\be\label{eq42}
&&\sum^{P}_{i=1}\Bigg\|\sum^{2}_{l=1}\left.\frac{\pa F_i}{\pa\g_l}
\right|_{(\g^{(m)}_0,\g^{(m)}_1,\la^{(m)}_1)}(\Delta\g^{(m)}_l)
+\left.\frac{\pa F_i}{\pa\la_1}\right|_{(\g^{(m)}_0,\g^{(m)}_1,\la^{(m)}_1)}(\Delta\la^{(m)}_1)
+F_i(\g^{(m)}_0,\g^{(m)}_1,\la^{(m)}_1)-u^\delta_{\infty,i}\Bigg\|^2_{L^2(S^1)}\no\\
&&\qquad\qquad\qquad+\beta\left(\sum^{2}_{l=1}||\Delta r^{(m)}_l||^2_{H^s(S^1)}
  +\sum^{2}_{l=1}|\Delta a^{(m)}_l|^2+(\Delta\la^{(m)}_1)^2\right)
\en
Here, $\beta\in\R^+$ is chosen such that
\be\label{eq43}
&&\left(\sum^{P}_{i=1}\Bigg\|\sum^{2}_{l=1}\left.\frac{\pa F_i}{\pa\g_l}
\right|_{(\g^{(m)}_0,\g^{(m)}_1,\la^{(m)}_1)}(\Delta\g^{(m)}_l)
+\left.\frac{\pa F_i}{\pa\la_1}\right|_{(\g^{(m)}_0,\g^{(m)}_1,\la^{(m)}_1)}(\Delta\la^{(m)}_1)
+F_i(\g^{(m)}_0,\g^{(m)}_1,\la^{(m)}_1)-u^\delta_{\infty,i}\Bigg\|_{L^2(S^1)}^2\right)^\frac12\no\\
&&\qquad\qquad\qquad=\rho\left(\sum^{P}_{i=1}\left\|F_i(\g^{(m)}_0,\g^{(m)}_1,\la^{(m)}_1)
-u^\delta_{\infty,i}\right\|_{L^2(S^1)}^2\right)^\frac12
\en
with a given constant $\rho<1$.
After obtaining $(\Delta\g^{(m)}_0,\Delta\g^{(m)}_1,\Delta\la^{(m)}_1)$,
the new approximation to $(\g_0,\g_1,\la_1)$ will be
$(\g^{(m+1)}_0,\g^{(m+1)}_1,\la^{(m+1)}_1)=(\g^{(m)}_0+\Delta\g^{(m)}_0,\g^{(m)}_1
+\Delta\g^{(m)}_1,\la^{(m)}_1+\Delta\la^{(m)}_1)$.
In the numerical computation, we can use the bisection algorithm to determine $\beta$
(see \cite{Hohage99}). Note that $\beta$ is unique if certain assumptions are satisfied
(see \cite{Hanke1997}). In this paper we choose $\beta$ by a posterior parameter choice rule.
For other strategies of solving (\ref{eq35}) we refer to \cite{Hohage99}.

Further, we define the relative error by
\be\label{eq44}
Err^{(m)}_{k_0}=\frac{1}{P}\sum^{P}_{i=1}\frac{\left\|F_i(\g^{(m)}_0,\g^{(m)}_1,\la^{(m)}_1)
-u^\delta_{\infty,i}\right\|_{L^2(S^1)}}{\left\|u^\delta_{\infty,i}\right\|_{L^2(S^1)}}
\en
Then the iteration is stopped if $Err^{(m)}_{k_0}<\tau\delta$,
where $\tau>1$ is a given constant.
Here, we use the subscript to emphasize its dependance on the wave number $k_0$.


In the numerical experiments, for fixed $d,k_0$ we use $u^\delta_\infty(\hat{x_i},d,k_0),i=1,\ldots,n,$
to denote the measured far field data, where $\hat{x_i}=2\pi(i-1)/n$.
We then seek the approximate solution $\g^{(m)}_0,\g^{(m)}_1,\la^{(m)}_1$ with
$r^{(m)}_0, r^{(m)}_1\in R_M\subset H^s(0,2\pi)$ ($s\ge0$) and $\la^{(m)}_1\in\R$,
where
\ben
R_M:=\big\{r\in H^s(0,2\pi)\,|\,r(\theta)=\al_0+\sum^{M}_{l=1}\Big[\al_l\cos(l\theta)
+\al_{l+M}\sin(l\theta)\Big],\,a_l\in\R,l=0,...,2M\big\}.
\enn
Similarly as in \cite{Hohage99}, we use the following:
\ben
\|f\|^2_{L^2(\Sp^1)}&\approx&\frac{2\pi}{n}\sum^{n}_{i=1}|f(\hat{x}_i)|^2,\quad f\in C(\Sp^1)\\
\|g\|_{H^s(\Sp^1)}^2&:=&2\pi\al_0+\pi\sum^{M}_{l=1}\Big[(1+l^2)^s(\al^2_l+\al^2_{l+M})\Big]
\enn
where $g(\theta)=\al_0+\sum^M_{l=1}[\al_l\cos(l\theta)+\al_{l+M}\sin(l\theta)].$
Then (\ref{eq42}), (\ref{eq43}), (\ref{eq44}) can be approximated, and the updated functions
are computed by determining their coefficients. For the determination of the coefficients
the reader is referred to \cite{Hohage99}, and for the computation of the Frechet derivative
of the far field operator we refer to \cite{HS,Hohage1998,Kirsch93}.

\begin{remark}\label{re3} {\rm
In Section \ref{sec2}, we observed that if $\la_1\to0$ then the solution $u$ of the
transmission problem (\ref{eq20})-(\ref{eq24}) and (\ref{rc}) converges to the solution
of the scattering problem (\ref{eq16})-(\ref{rc}) with a Neumann boundary condition
$\mathscr{B}(u):=\pa u/\pa\nu=0$ on $S_1$ and that if $\la_1\to\infty$ (or $\tau_1\to0$)
then the solution $u$ of the transmission problem (\ref{eq20})-(\ref{eq24}) and (\ref{rc})
converges to the solution of the scattering problem (\ref{eq16})-(\ref{rc}) with a Dirichlet
boundary condition $\mathscr{B}(u):=u=0$ on $S_1$.
Thus, the following fact is expected: 1) when $u_\infty$ is the far field pattern corresponding to a
sound-hard obstacle $\Om_2$ (i.e., a Neumann boundary condition on $S_1$) then $\lambda^{(m)}_1$
given by the Newton iteration method for (\ref{eq35}) will be very close to zero, and
2) when $u_\infty$ is the far field pattern corresponding to a sound-soft obstacle $\Om_2$
(i.e., a Dirichlet boundary condition on $S_1$) then $\tau^{(m)}_1$ given by the Newton iteration
method for (\ref{eq36}) will be very close to zero.
However, there is no rigorous proof of this convergence result even for the case of
an obstacle in the homogeneous background medium.
}
\end{remark}

\begin{remark}\label{re4} {\rm
Since we do not know the boundary condition on $S_1$, we choose a constant $\tilde{\la}$. 
Let $\lambda^{(m)}_1$ be the approximation of $\lambda_1$ after the $m$-th iteration and
$\tau^{(m)}_1=1/\lambda^{(m)}_1$. Then, if $|\la^{(m)}_1|<=\tilde{\la}$, we solve (\ref{eq35})
to get $(\g^{(m+1)}_0,\g^{(m+1)}_1,\la^{(m+1)}_1)$ and let $\tau^{(m+1)}_1=1/\la^{(m+1)}_1$. 
On the other hand, if $|\la^{(m)}_1|>\tilde{\la}$ then we solve (\ref{eq36}) to get 
$(\g^{(m+1)}_0,\g^{(m+1)}_1,\tau^{(m+1)}_1)$ and let $\la^{(m+1)}_1=1/\tau^{(m+1)}_1$. 
It should be noted that, in numerical computation, $\la^{(m)}_1$ may be negative. 
This does not affect the iteration process, and $\la^{(m)}_1$ can still be used to
determine the type of boundary conditions on $S_1$. 
}
\end{remark}

\begin{remark}\label{re5} {\rm
In our iteration algorithm, multiple frequency data are considered in order to get better 
reconstructions. 
General speaking, the low frequency data are used to get the main profiles of
the boundaries and the higher frequency data are used to get a refined reconstruction.
}
\end{remark}

Our Newton iteration algorithm is given in Algorithm \ref{alg1} below.

\begin{algorithm}\label{alg1}
Given $\la_0,n_1,u^\delta_\infty(\cdot,d_i,k_{0j}),i=1,...,P,j=1,...,Q$,
where $k_{01}<k_{02}<\cdots<k_{0P}.$
\begin{description}
\item 1) Choose initial guesses $\g^{(0)}_0,\g^{(0)}_1,\la^{(0)}_1$ for
$\g_0,\g_1,\la_1$, where $\la^{(0)}_1\neq 0$, and let $Ir$ be the maximal iteration number 
for each fixed frequency. Set $j=0$, $\tau^{(0)}_1=1/\la^{(0)}_1$ and go to Step 2).
\item 2) Set $j=j+1$. If $j>Q$, then stop the iteration; otherwise, set $k_0=k_{0j},m=0$ 
and go to Step 3).
\item 3) If $Err^{(m)}_{k_0}<\tau\delta$ or $m=Ir$, then set
$(\g^{(0)}_0,\g^{(0)}_1,\la^{(0)}_1,\tau^{(0)}_1)=(\g^{(m)}_0,\g^{(m)}_1,\la^{(m)}_1,\tau^{(m)}_1)$ 
and return to Step 2; otherwise, go to Step 4).
\item 4) Using the strategy in Remark \ref{re4}, we get the new approximate data
$(\g^{(m+1)}_0,\g^{(m+1)}_1,\la^{(m+1)}_1,\tau^{(m+1)}_1)$.
\item 5) Set $m=m+1$, and return to Step 3).
\end{description}
\end{algorithm}


\begin{remark}\label{re7} {\rm
Let $\la_1^{(end)}$ be the final approximation of $\la_1$. If
$|\la^{(end)}_1|>> 1$, then it is concluded that $S_1$ is a sound-soft obstacle,
and if $|\la^{(end)}_1|<< 1$, we then conclude that $S_1$ is a sound-hard obstacle.
}
\end{remark}


\section{Numerical examples}\label{sec5}
\setcounter{equation}{0}

In this section, we present several numerical experiments to demonstrate the effectiveness 
of our algorithm. In all numerical experiments we make the following assumptions.
\begin{description}
\item 1) We choose $n_1=0.64,\la_0=1.2,s=1.6$ (which is suggested in \cite{Hohage99}),
M=25, $\rho=0.8$, $\tau=1.5$ and $\tilde{\la}=1$ for iteration.
Since $k_2$ is a parameter at our disposal, we just choose $k_2=k_1$ for each frequency. 
\item 2) To generate the synthetic data, we use the integral equation method in Section \ref{sec3}
to solve the direct problem. 
In order to avoid inverse crimes, we choose $n=128$ for generating the synthetic data and $n=64$
for solving the direct problem at each iteration.
\item 3) We use the full far-field data with $64$ measurement points.
The noisy data $u^\delta_\infty$ is obtained as 
$u^\delta_\infty=u_\infty+\delta\zeta{||u_\infty||_{L^2(\Sp^1)}}/{||\zeta||_{L^2(\Sp^1)}}$,
where $\zeta$ is a random number with $\text{Re}(\zeta),\text{Im}(\zeta)\in N(0,1)$.
\item 4) The initial shapes of $S_0,S_1$ are chosen to be circles.
\item 5) For the incident plane waves, we choose $d_i=(cos(\theta_i),sin(\theta_i)),i=1,...,P,$ 
where $P$ is the number of incident plane waves and $\theta_i=2\pi(i-1)/P$.
\item 6) In each figure, we use solid line '-', dotted line '...' and dashed line '- -' 
to represent the actual curve, the initial guess of the curve and the reconstructed curve, 
respectively.
\item 7) The shapes of the actual boundaries are given in Table \ref{table2}.
\end{description}

\begin{table}[h]
\centering
\begin{tabular}{ll}
\hline
Type & Parametrization\\
\hline
Circle & $r_0(\cos{t},\sin{t}),\;t\in[0,2\pi]$\\
Apple shaped & $[({0.5+0.4\cos{t}+0.1\sin(2t)})/({1+0.7\cos{t}})](\cos{t},\sin{t}),\;t\in[0,2\pi]$\\
Kite shaped & $(\cos{t}+0.65\cos(2t)-0.65,1.5\sin{t}),\;t\in[0,2\pi]$\\
Rounded square & $({3}/{2})(\cos^3{t}+\cos{t},\sin^3{t}+\sin{t}),\;t\in[0,2\pi]$\\
Rounded triangle & $(2+0.3\cos(3t))(\cos{t},\sin{t}),\;t\in[0,2\pi]$\\
\hline
\end{tabular}
\caption{Parametrization of the boundaries $S_0$ and $S_1$.}\label{table2}
\end{table}

\begin{figure}[htbp]
\centering
\subfigure[\textbf{no noise}, $\la^{(Ir)}_1=-502.6$, $Err^{(Ir)}=2.151\times10^{-3}$, $Ir=25$]
{\label{fig1:subfig:a}
\includegraphics[width=3in]{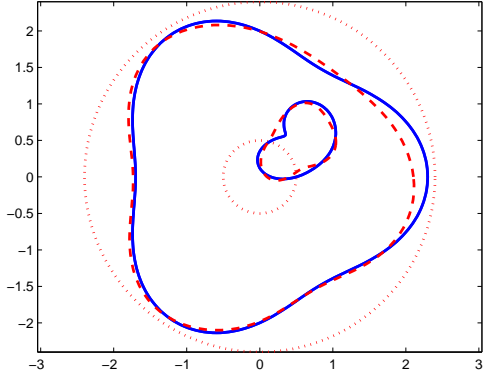}}
\subfigure[\textbf{no noise}, $\la^{(Ir)}_1= 1004$, $Err^{(Ir)}=2.246\times10^{-3}$, $Ir=25$]
{\label{fig1:subfig:b}
\includegraphics[width=3in]{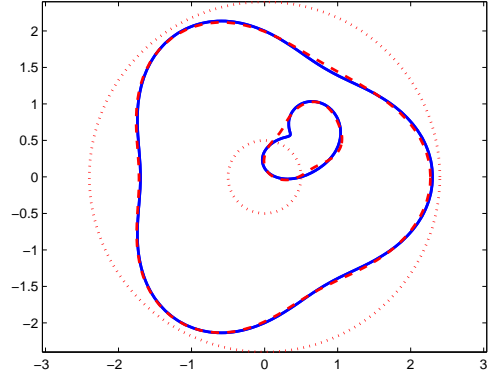}}
\subfigure[\textbf{no noise}, $\la^{(Ir)}_1=-9971$, $Err^{(Ir)}=1.861\times10^{-3}$, $Ir=25$]
{\label{fig1:subfig:c}
\includegraphics[width=3in]{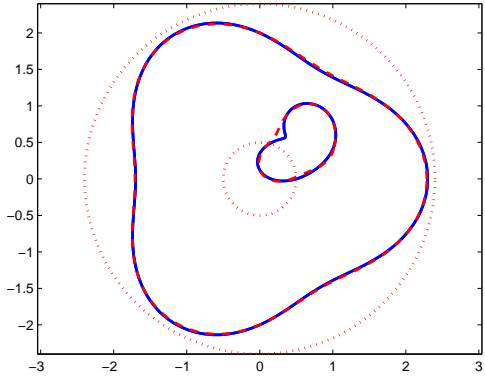}}
\subfigure[\textbf{no noise}, $\la^{(Ir)}_1= 2796$, $Err^{(Ir)}=1.975\times10^{-3}$, $Ir=25$]
{\label{fig1:subfig:d}
\includegraphics[width=3in]{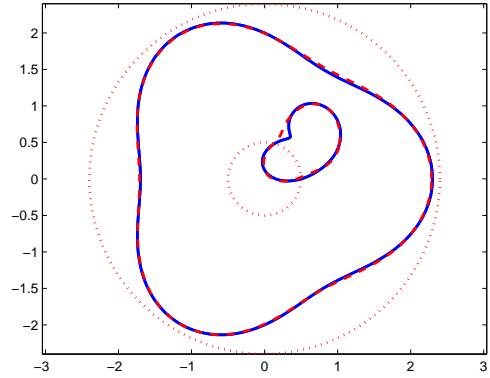}}
\caption{An apple-shaped, sound-soft obstacle embedded in a rounded triangle-shaped penetrable 
obstacle. The boundary condition on $S_1$ is assumed unknown, and the exact data with a fixed 
wave number $k_0=2$ are used. (a), (b), (c) and (d) give the reconstructed results with using 
one, two, three and four incident plane waves, respectively.}\label{fig1:subfig} 
\end{figure}

\begin{figure}[htbp]
\centering
\subfigure[\textbf{3\% noise}, $k_0=2$, $\la^{(Ir)}_1= -275.4$, $Err^{(Ir)}=0.04053$, $Ir=12$]
{\label{fig1b:subfig:a} 
\includegraphics[width=3in]{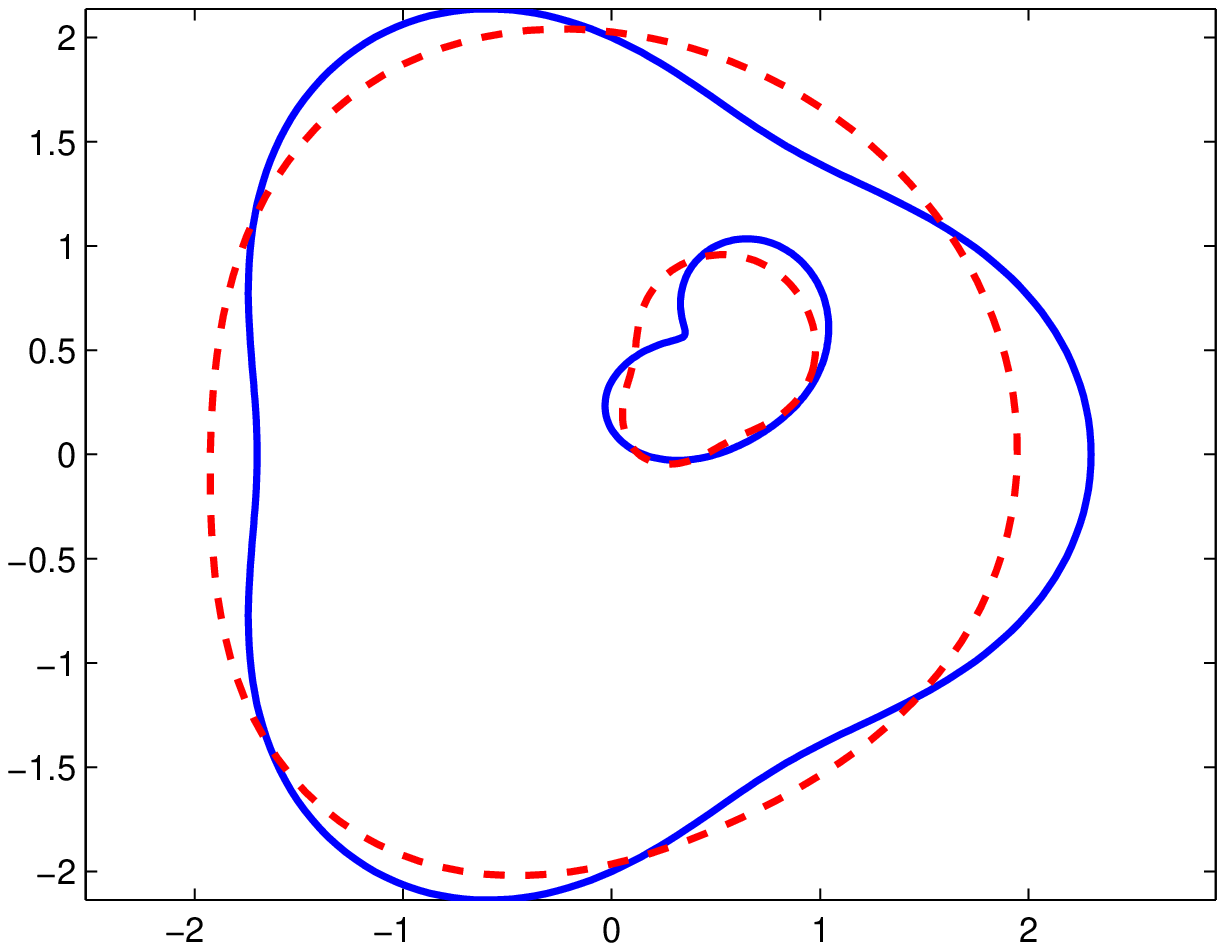}}
\hspace{0in}
\subfigure[\textbf{3\% noise}, $k_0=4$, $\la^{(Ir)}_1=-515$, $Err^{(Ir)}=0.03803$, $Ir=6$]
{\label{fig1b:subfig:b} 
\includegraphics[width=3in]{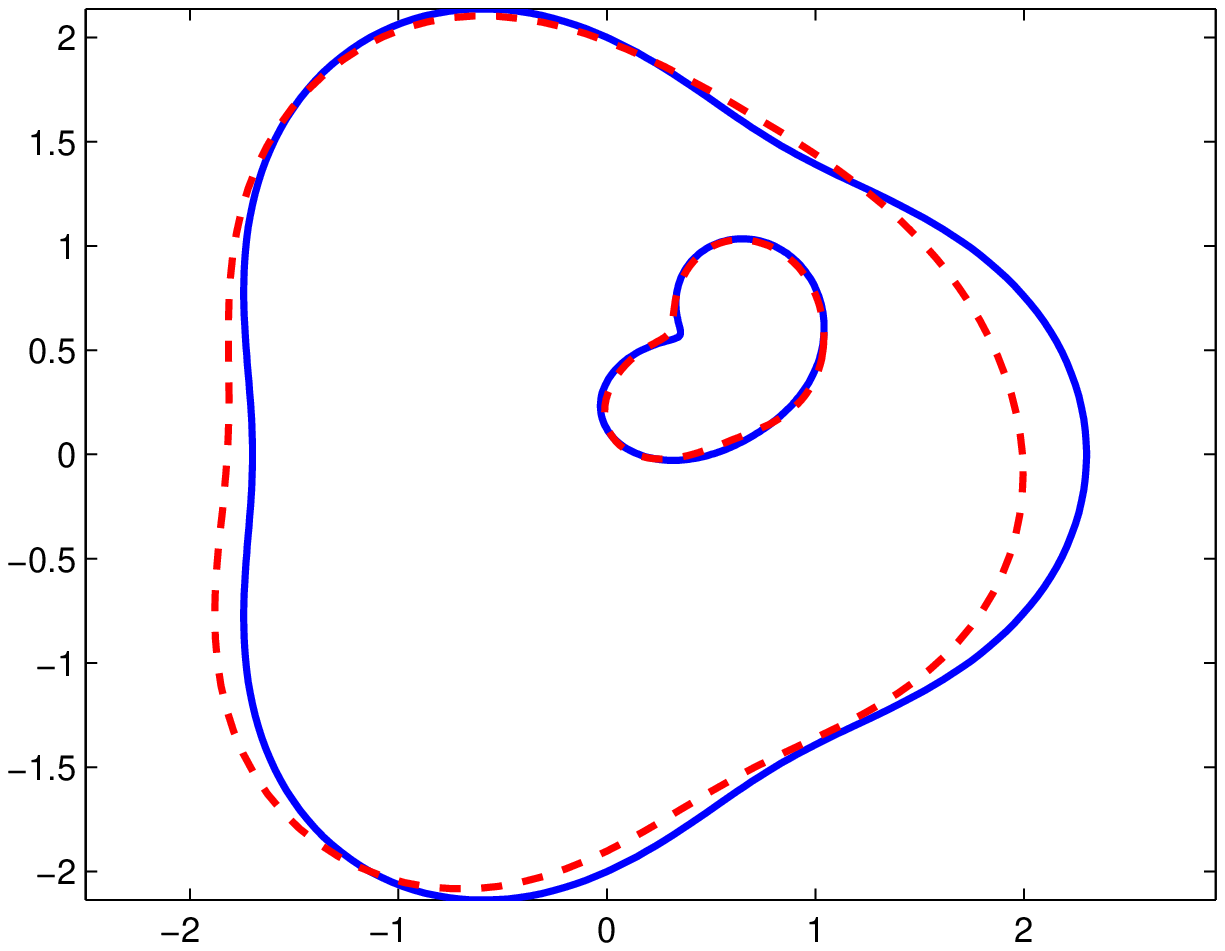}}
\subfigure[\textbf{3\% noise}, $k_0=6$, $\la^{(Ir)}_1=736.2$, $Err^{(Ir)}=0.04316$, $Ir=2$]
{\label{fig1b:subfig:c} 
\includegraphics[width=3in]{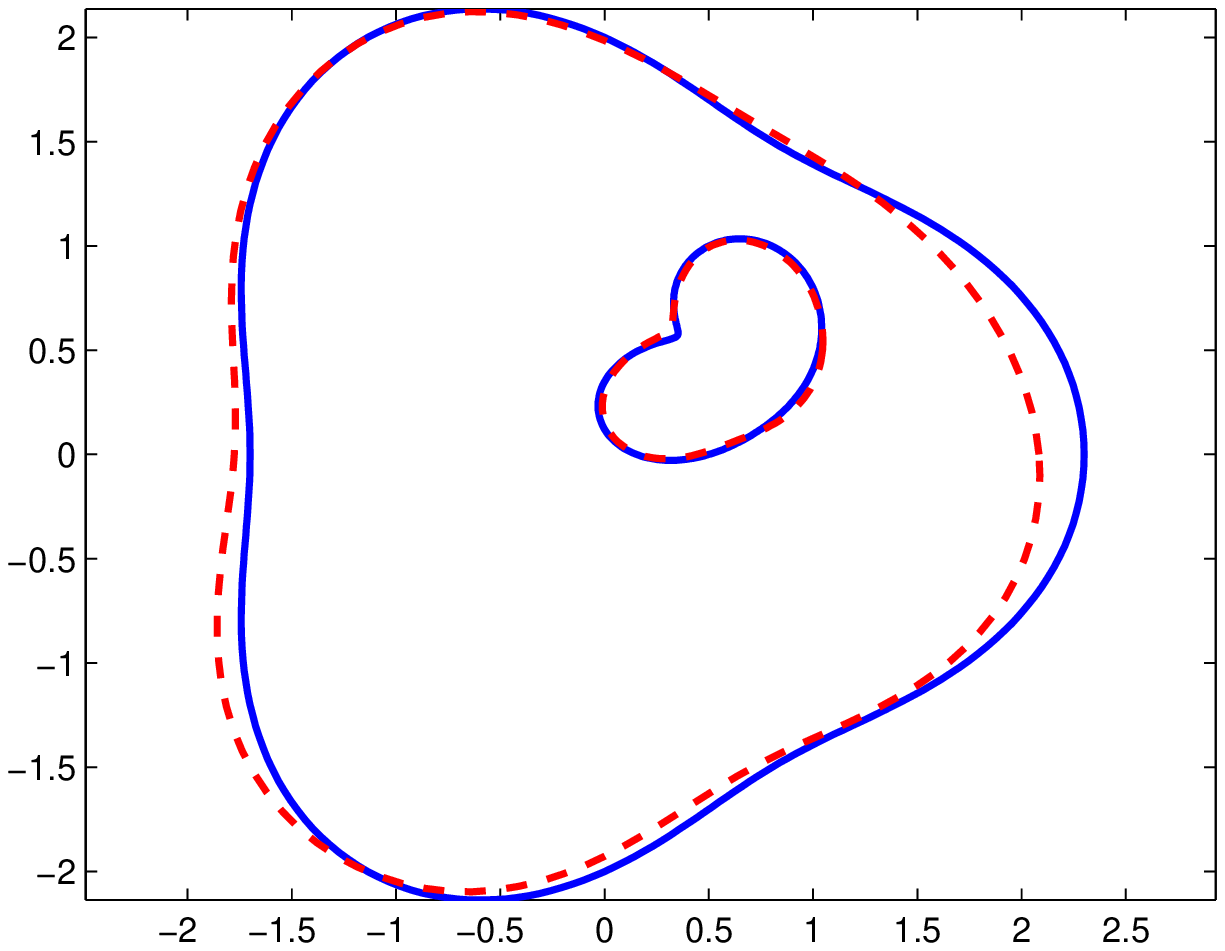}}
\hspace{0in}
\subfigure[\textbf{3\% noise}, $k_0=8$, $\la^{(Ir)}_1=-895.7$, $Err^{(Ir)}=0.03651$, $Ir=4$]
{\label{fig1b:subfig:d} 
\includegraphics[width=3in]{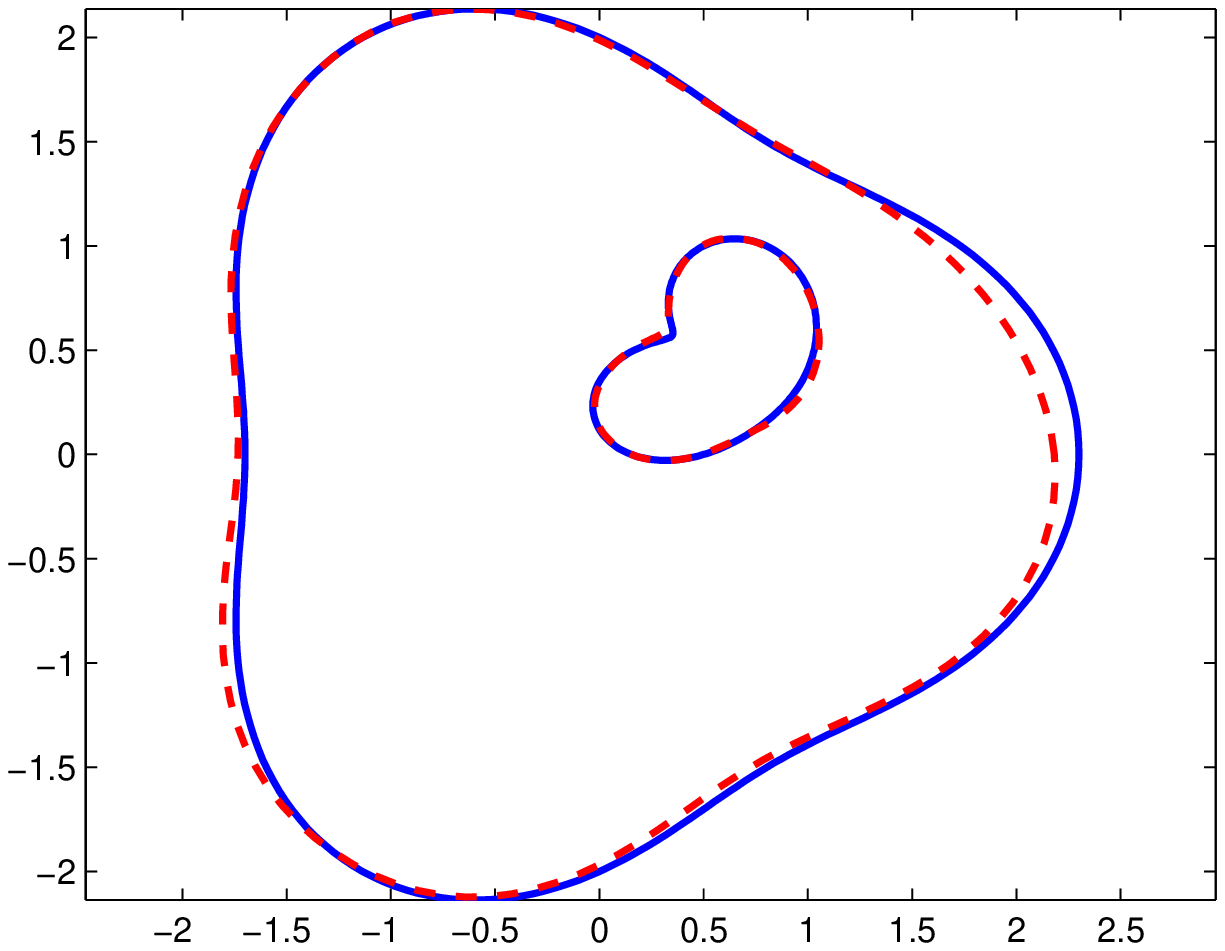}}
\caption{An apple-shaped, sound-soft obstacle embedded in a rounded triangle-shaped penetrable 
obstacle. The boundary condition on $S_1$ is assumed unknown, and noisy multi-frequency data with 
one incident plane wave are used. (a), (b), (c) and (d) show the reconstructions at $k_0=2,4,6,8,$ 
respectively.}\label{fig1b:subfig} 
\end{figure}

\begin{figure}[htbp]
\centering
\subfigure[\textbf{3\% noise}, $k_0=2$, $\la^{(Ir)}_1=-2591$, $Err^{(Ir)}=0.03686$, $Ir=12$]
{\label{fig1c:subfig:a} 
\includegraphics[width=3in]{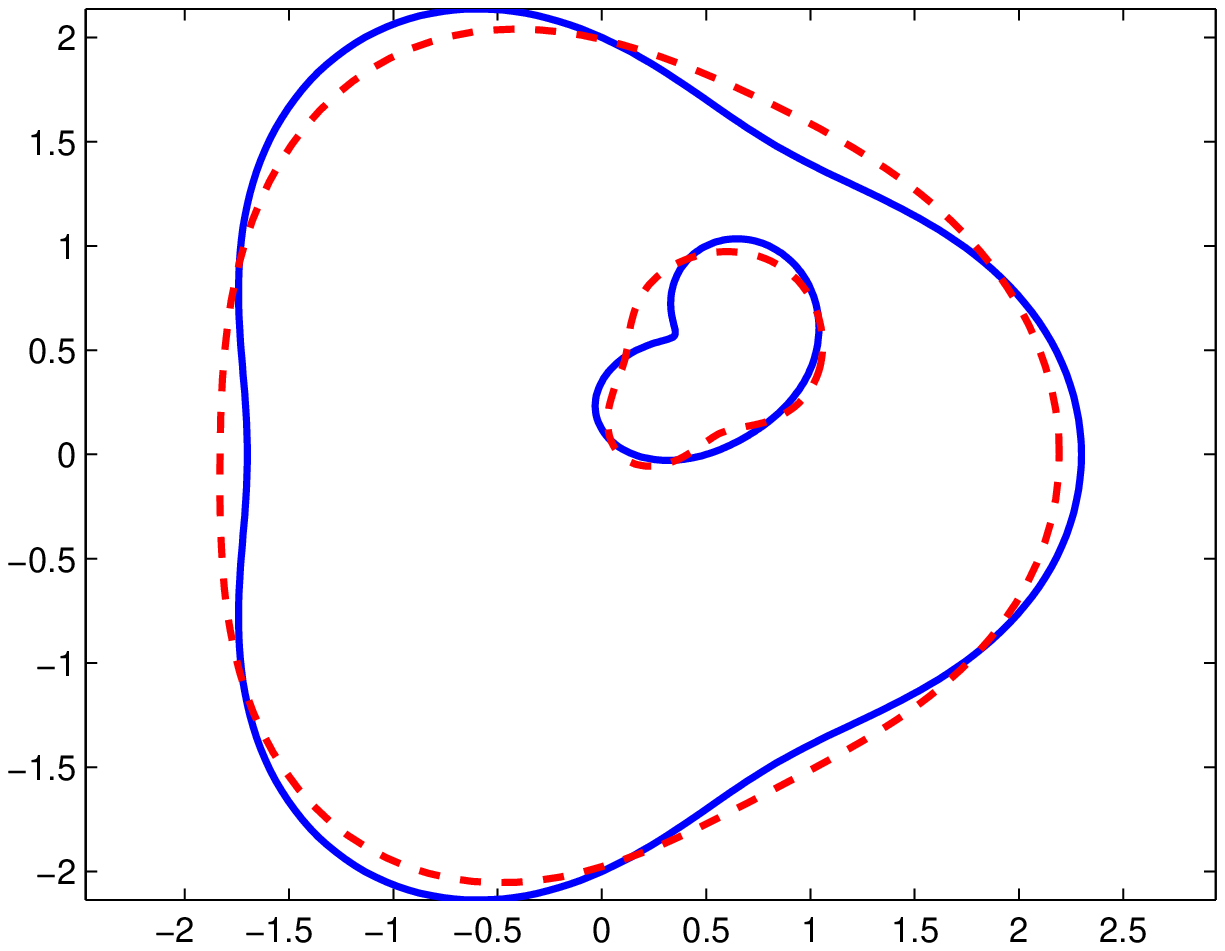}}
\hspace{0in}
\subfigure[\textbf{3\% noise}, $k_0=4$, $\la^{(Ir)}_1=963.5$, $Err^{(Ir)}=0.03987$, $Ir=4$]
{\label{fig1c:subfig:b} 
\includegraphics[width=3in]{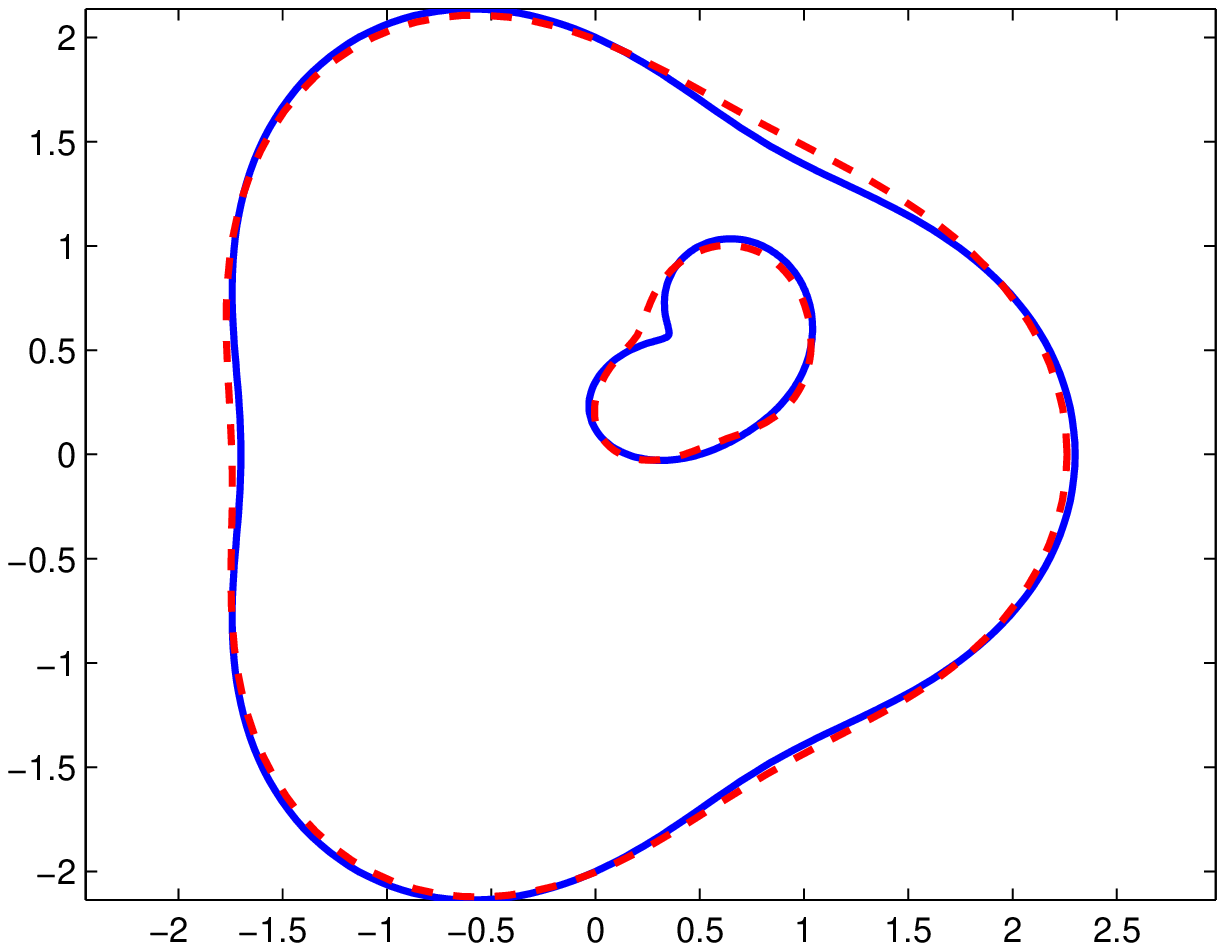}}
\subfigure[\textbf{3\% noise}, $k_0=6$, $\la^{(Ir)}_1= -417.3$, $Err^{(Ir)}=0.04476$, $Ir=2$]
{\label{fig1c:subfig:c} 
\includegraphics[width=3in]{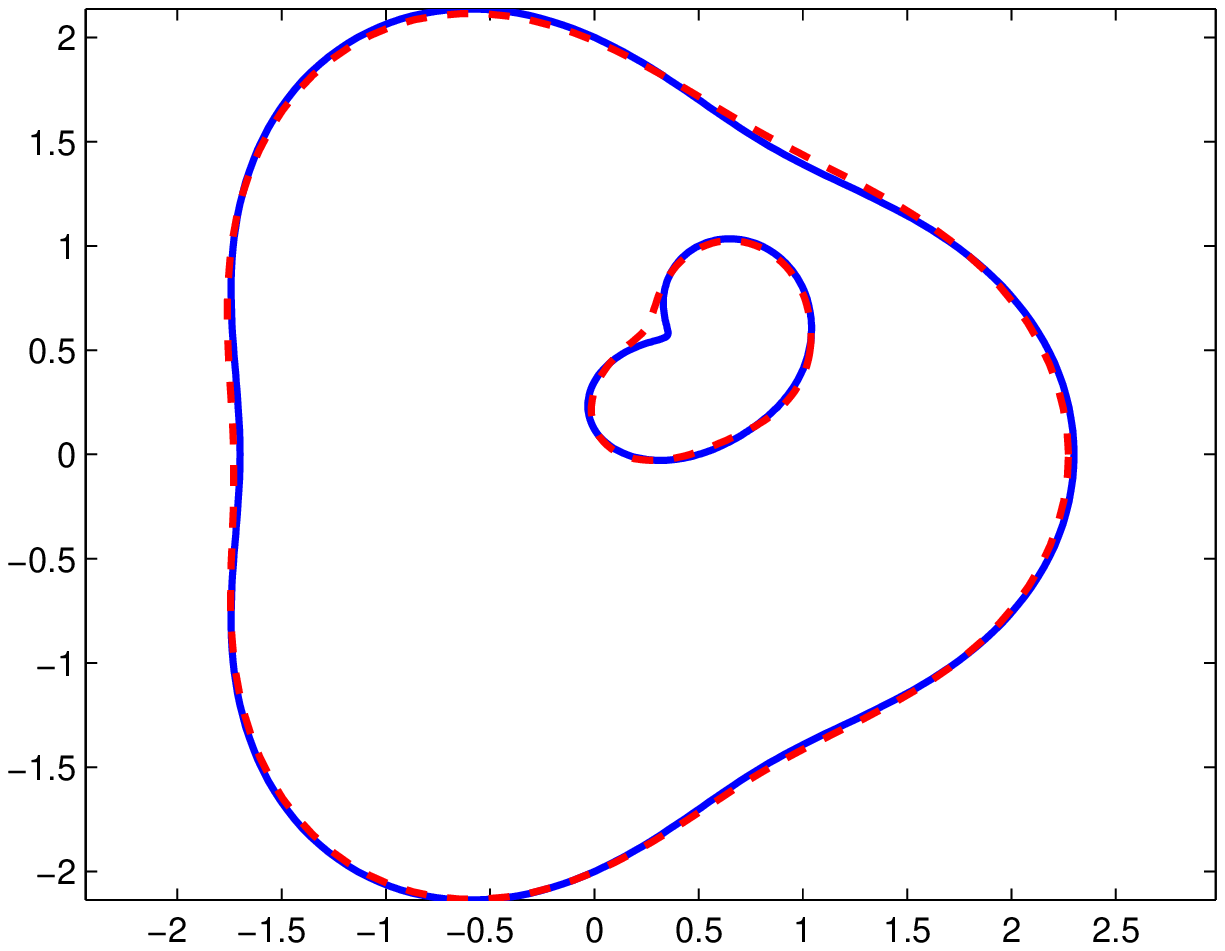}}
\hspace{0in}
\subfigure[\textbf{3\% noise}, $k_0=8$, $\la^{(Ir)}_1=5715$, $Err^{(Ir)}=0.04139$, $Ir=2$]
{\label{fig1c:subfig:d} 
\includegraphics[width=3in]{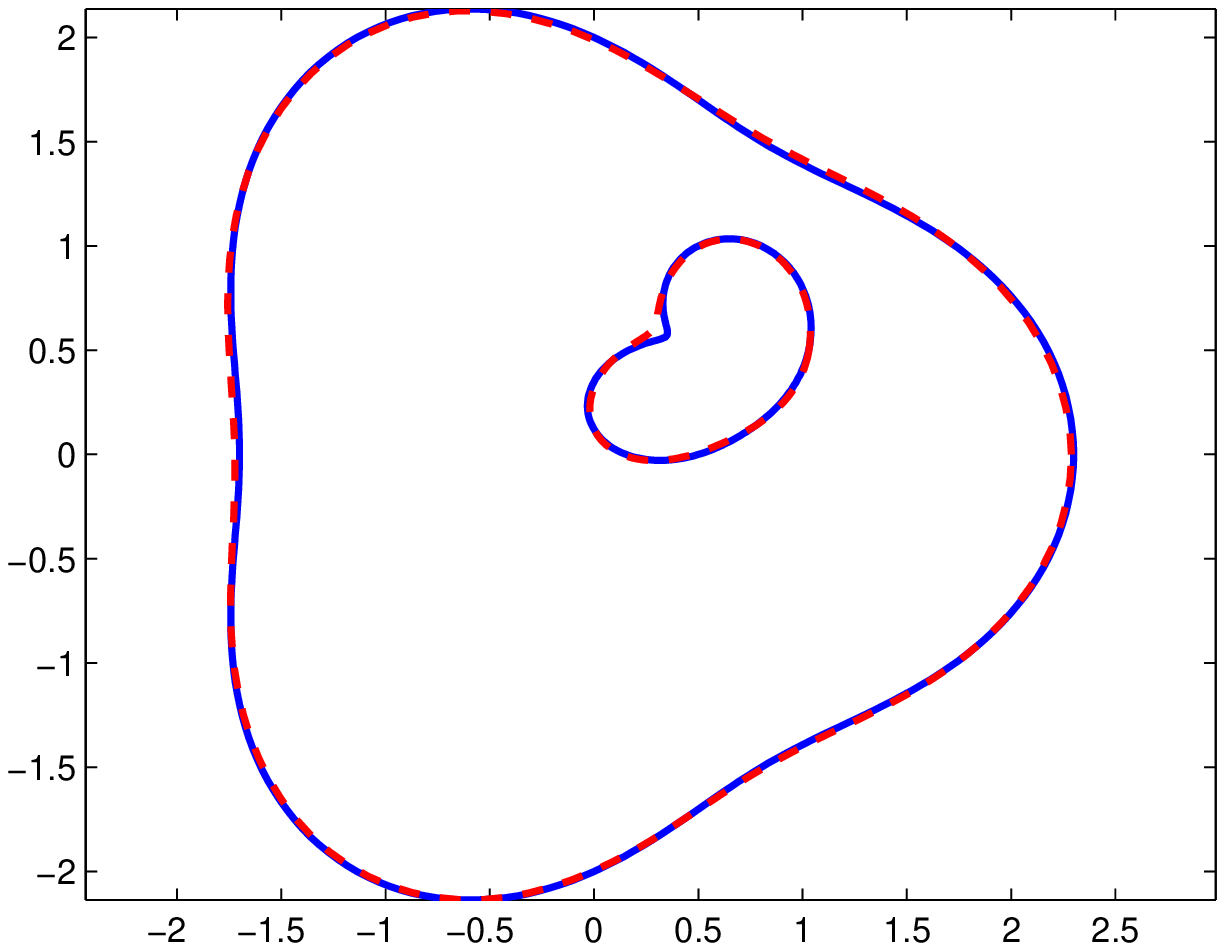}}
\caption{An apple-shaped, sound-soft obstacle embedded in a rounded triangle-shaped penetrable 
obstacle. The boundary condition on $S_1$ is assumed unknown, and noisy multi-frequency data with 
four incident plane waves are used. (a), (b), (c) and (d) show the reconstructions at $k_0=2,4,6,8,$ 
respectively.}\label{fig1c:subfig} 
\end{figure}

\textbf{Example 1.} In this example, we consider a sound-soft, apple-shaped obstacle embedded 
in a rounded triangle-shaped penetrable obstacle.
Suppose the far-field data are given for this case which corresponds to a Dirichlet boundary 
condition on the boundary $S_1$ of the embedded obstacle.
However, the boundary condition on $S_1$ is assumed unknown, and we will use our 
method to reconstruct both obstacles and determine the type of boundary conditions on $S_1$.
For the initial guesses $\g^{(0)}_0,\g^{(0)}_1$ of $S_0,S_1$, respectively, we choose
the radii $r^{(0)}_0=2.4,\;r^{(0)}_1=0.5$, respectively, and the two centers to be the origin 
$(0,0)$. By trials, choose $\la^{(0)}=10$. Here, the center of $\g^{(m)}_1$ is updated at 
each iteration. In Figure \ref{fig1:subfig}, the reconstruction results are given by using 
the exact data with a fixed wave number $k_0=2$. 
Figure \ref{fig1:subfig:a}, \ref{fig1:subfig:b}, \ref{fig1:subfig:c} and \ref{fig1:subfig:d} 
show the reconstructed results with one, two, three and four incident plane waves, respectively. 
It is observed that using more incident plane waves yields a better reconstruction result. 
Figures \ref{fig1b:subfig} and \ref{fig1c:subfig} give the reconstruction results by using
the $3\%$ noisy far-field data with one and four incident plane waves, respectively. Here, we
use multiple frequency data with $k_0=2,4,6,8$. From Figures \ref{fig1b:subfig} and 
\ref{fig1c:subfig} it is seen that multiple frequency data yield excellent reconstruction results. 
Finally, since $|\la^{(Ir)}_1|>>1$ for each case, we can conclude that $S_1$ is a sound-soft 
boundary.

\begin{figure}[htbp]
\centering
\subfigure[\textbf{3\% noise}, $Err^{(Ir)}=0.04301$, $Ir=12$]
{\label{fig2:subfig:a} 
\includegraphics[width=3in]{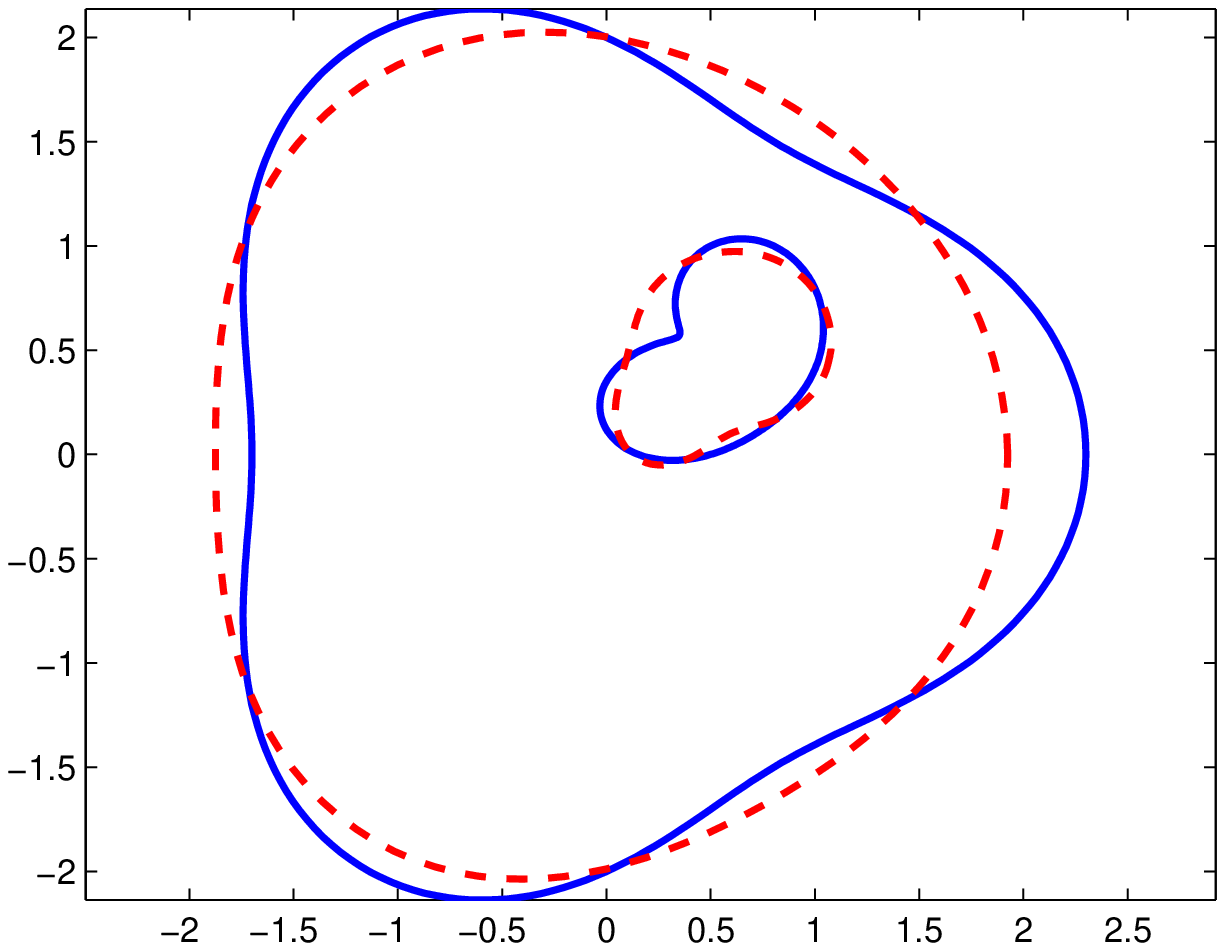}}
\hspace{0in}
\subfigure[\textbf{3\% noise}, $Err^{(Ir)}=0.04308$, $Ir=5$]
{\label{fig2:subfig:b} 
\includegraphics[width=3in]{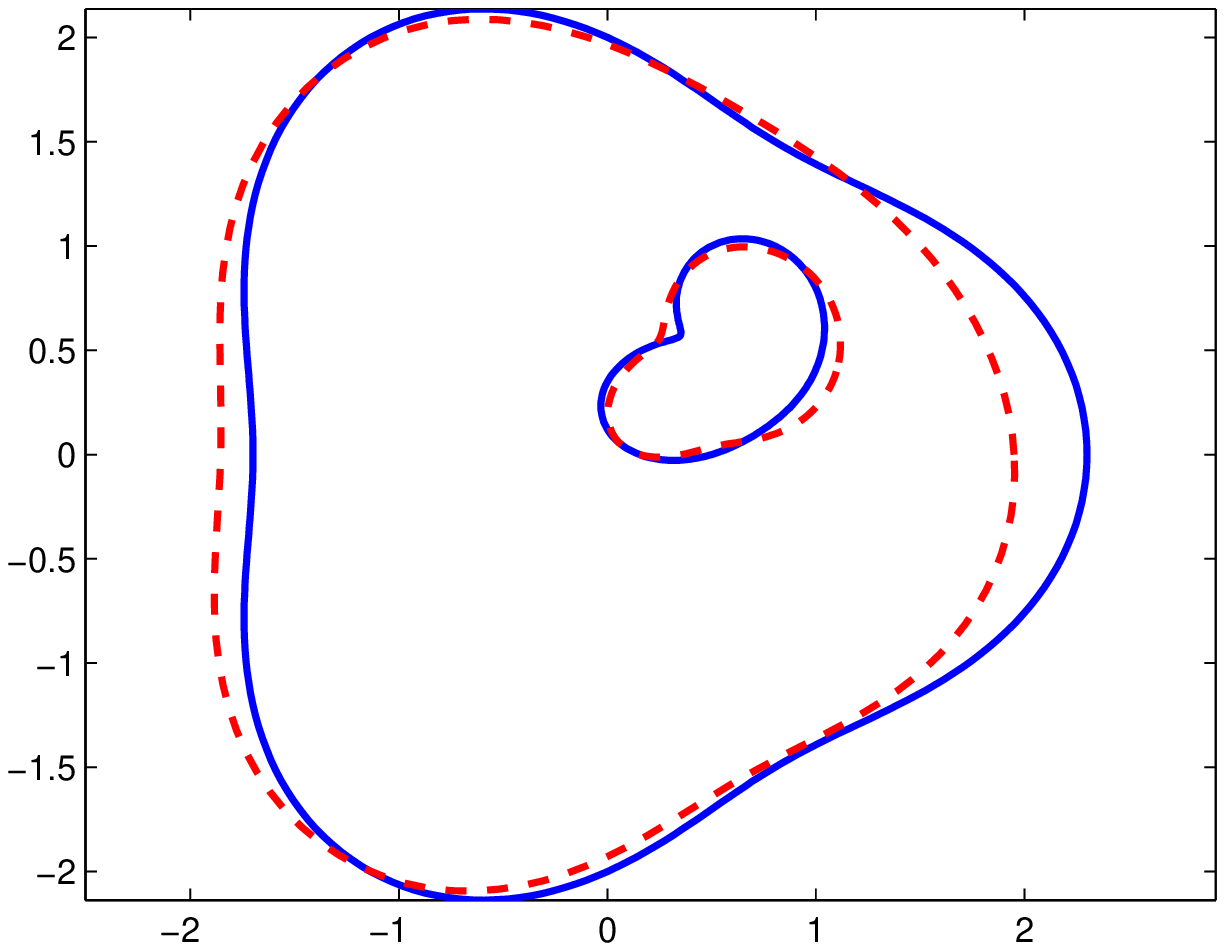}}
\subfigure[\textbf{3\% noise}, $Err^{(Ir)}=0.03717$, $Ir=3$]
{\label{fig2:subfig:c} 
\includegraphics[width=3in]{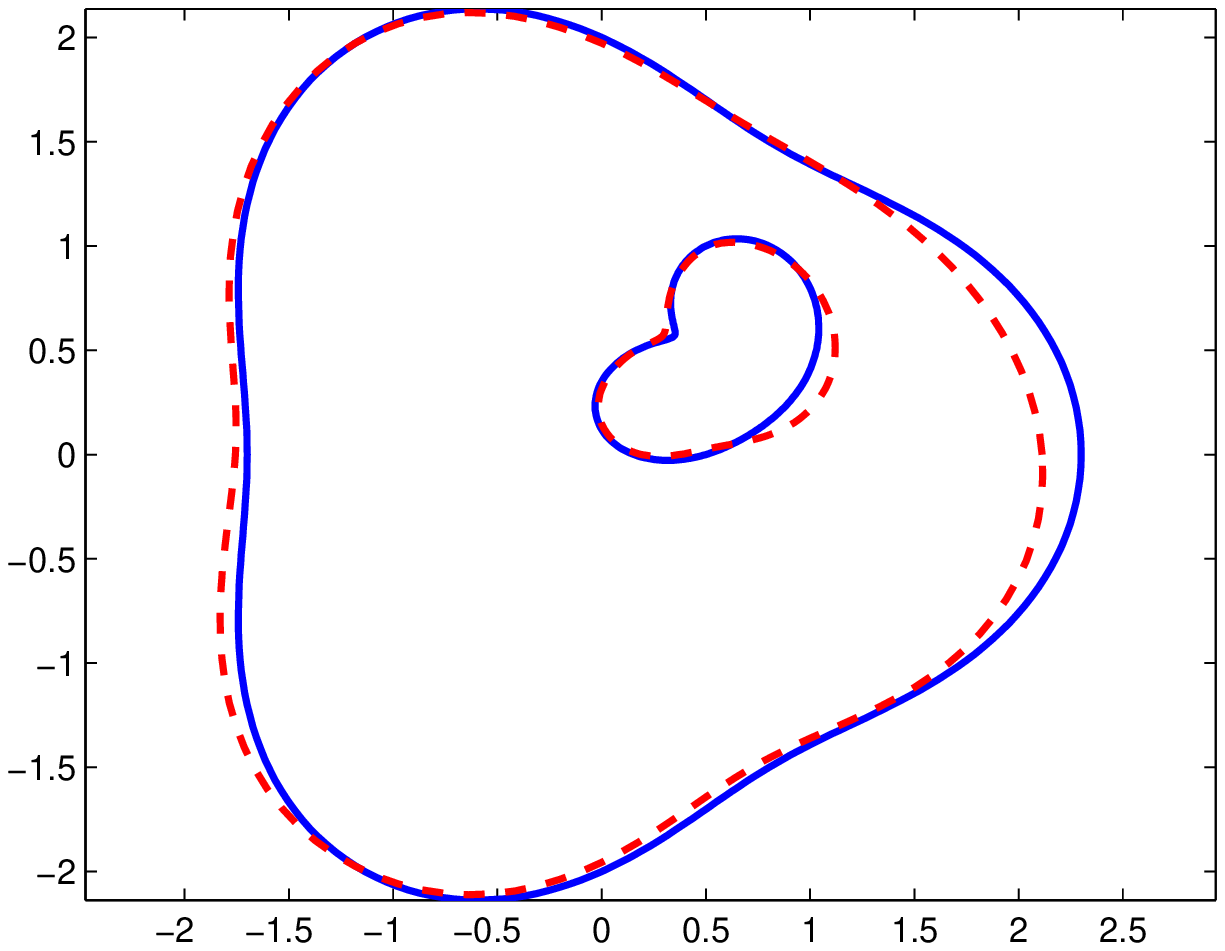}}
\hspace{0in}
\subfigure[\textbf{3\% noise}, $Err^{(Ir)}=0.03687$, $Ir=2$]
{\label{fig2:subfig:d} 
\includegraphics[width=3in]{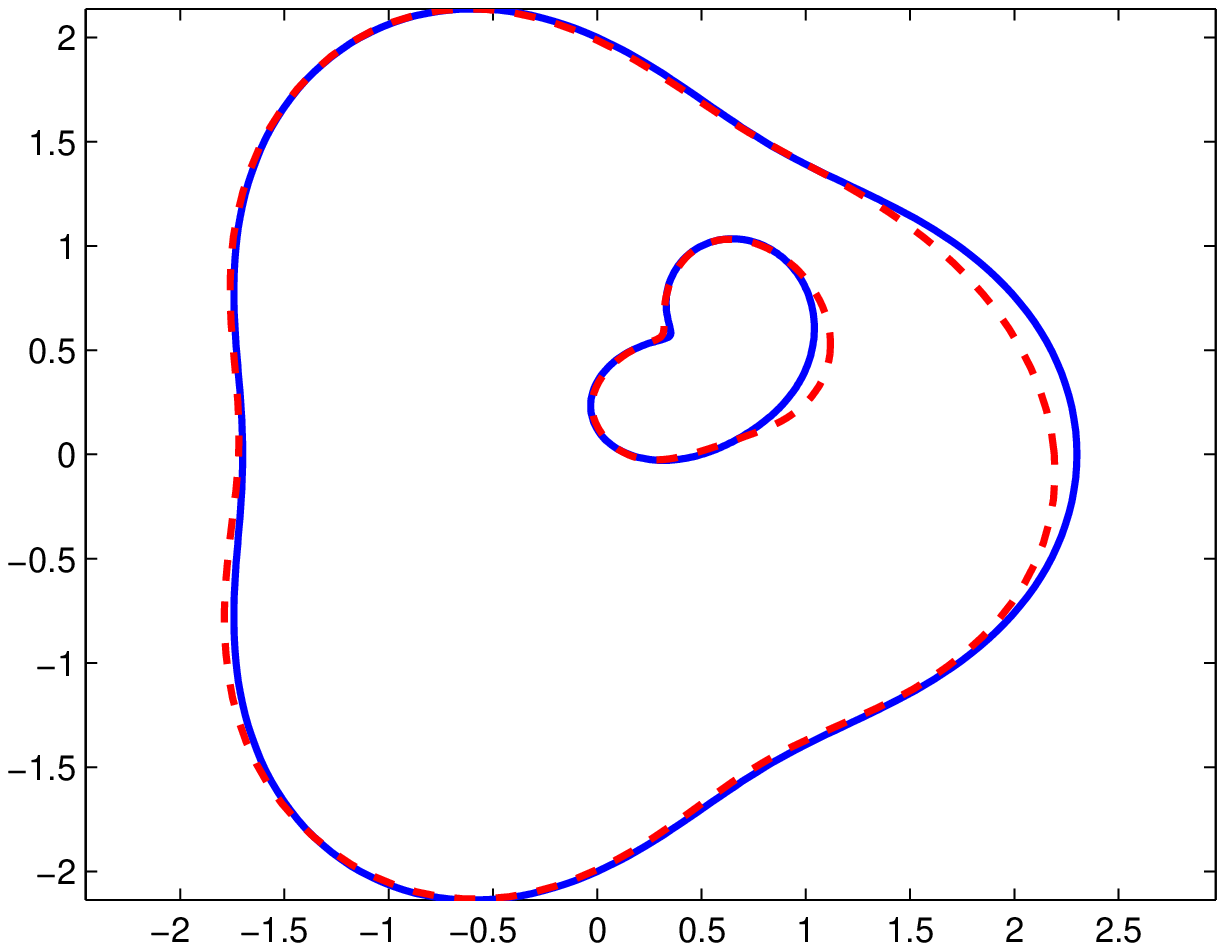}}
\caption{An apple-shaped, sound-soft obstacle embedded in a rounded triangle-shaped penetrable
obstacle. The boundary condition on $S_1$ is assumed known, and noisy multi-frequency data with 
one incident plane wave are used. (a), (b), (c) and (d) show the reconstructions at $k_0=2,4,6,8,$
respectively.}\label{fig2:subfig} 
\end{figure}

\textbf{Example 2.} In this example, we consider the problem in Example 1 again.
Here, the far-field data are given again for the case of an embedded, sound-soft obstacle,
and $S_1$ is also assumed to be a sound-soft boundary.
We will use the Newton iteration method for the inverse problem corresponding to the original 
scattering problem (\ref{eq16})-(\ref{rc}) with $\mathscr{B}(u):=u$ to reconstruct the two obstacles.
The initial guesses of $\g_0,\g_1,\la_0$ are the same as in Example 1, and the center of 
$\g^{(m)}_1$ is updated at each iteration.
In Figure \ref{fig2:subfig}, we present the reconstruction results by using
the $3\%$ noisy far-field data with one incident direction and multiple frequencies $k_0=2,4,6,8$. 
From Figures \ref{fig1b:subfig} and \ref{fig2:subfig} it can be seen that the reconstruction 
obtained by both methods is almost the same.

\begin{figure}[htbp]
\centering
\subfigure[\textbf{no noise}, $\la^{(Ir)}_1=4.186\times10^{-4}$, $Err^{(Ir)}=0.001296$, $Ir=30$]
{\label{fig3:subfig:a}
\includegraphics[width=3in]{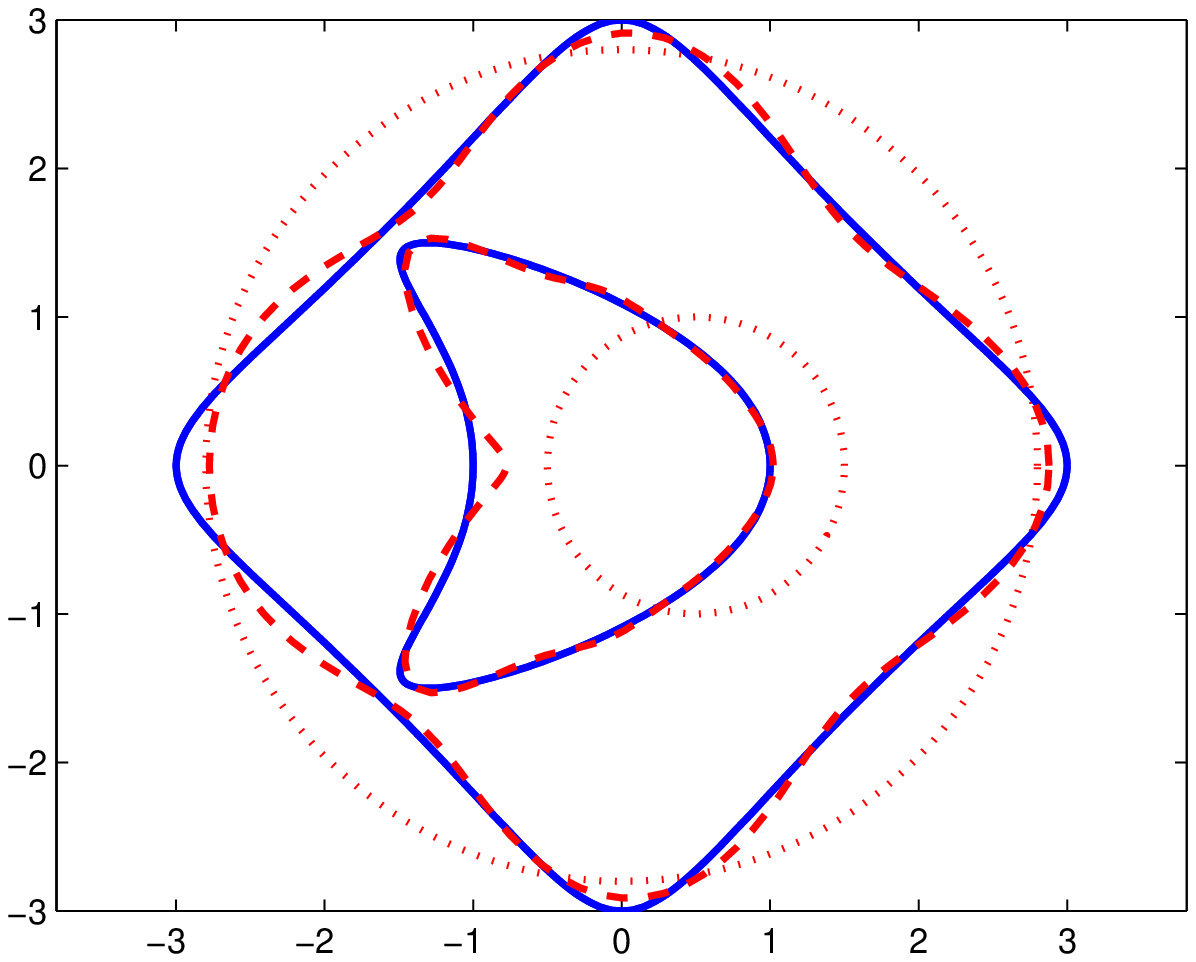}}
\hspace{0in}
\subfigure[\textbf{no noise}, $\la^{(Ir)}_1= 2.097\times10^{-5}$, $Err^{(Ir)}=0.01086$, $Ir=30$]
{\label{fig3:subfig:b}
\includegraphics[width=3in]{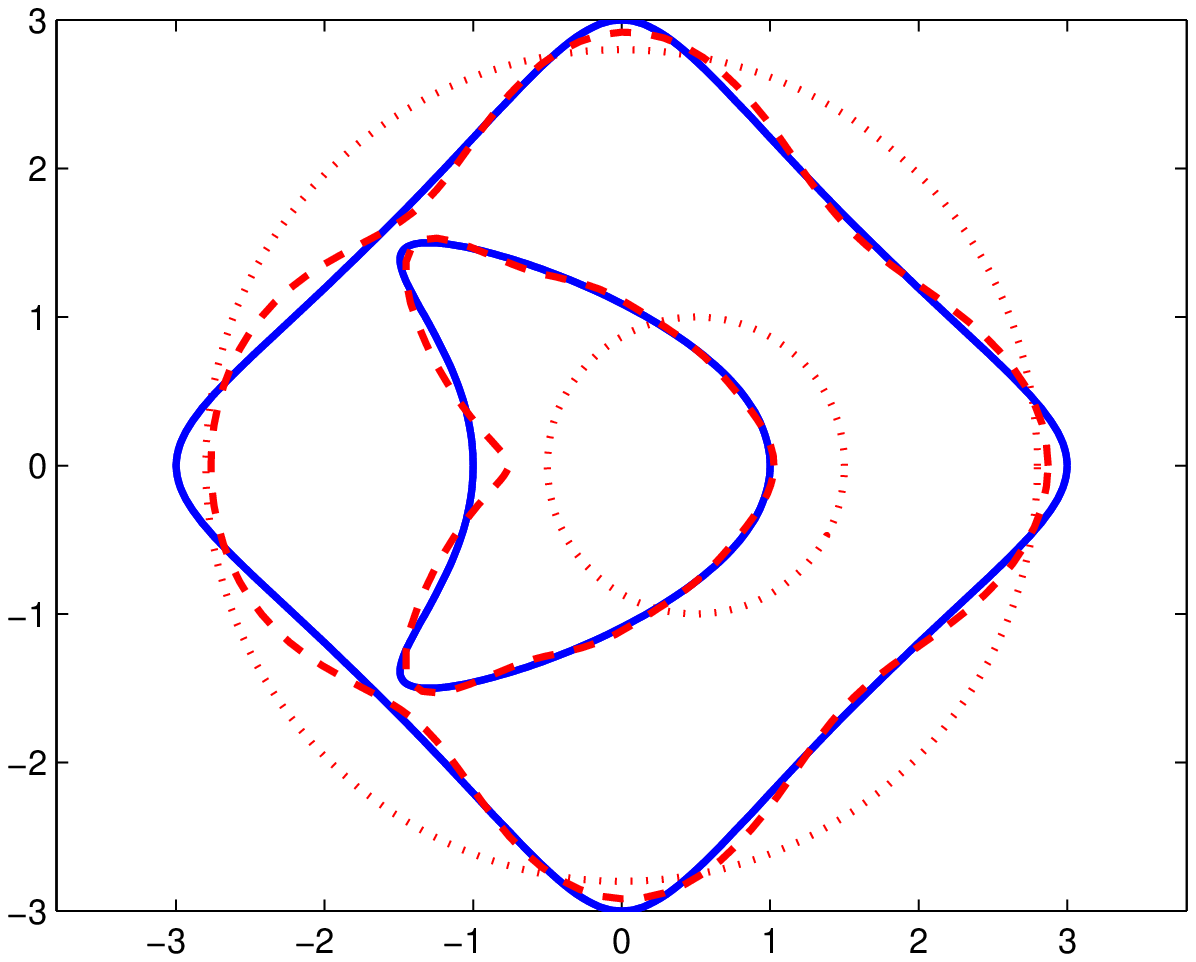}}
\subfigure[\textbf{no noise}, $\la^{(Ir)}_1=6.035\times10^{-4}$, $Err^{(Ir)}=0.001196$, $Ir=30$]
{\label{fig3:subfig:c} 
\includegraphics[width=3in]{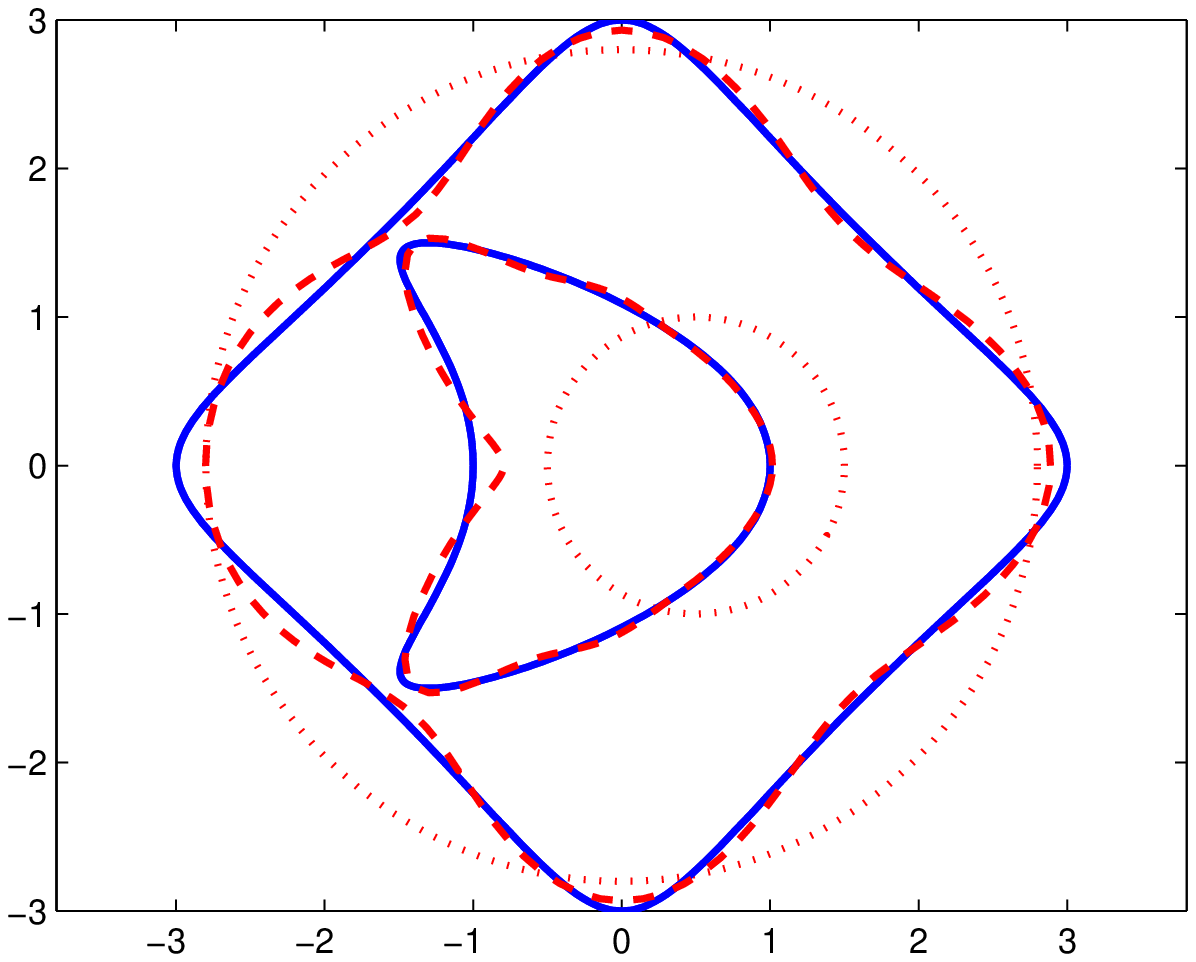}}
\hspace{0in}
\subfigure[\textbf{no noise}, $\la^{(Ir)}_1=4.963\times10^{-4}$, $Err^{(Ir)}=0.001212$, $Ir=30$]
{\label{fig3:subfig:d} 
\includegraphics[width=3in]{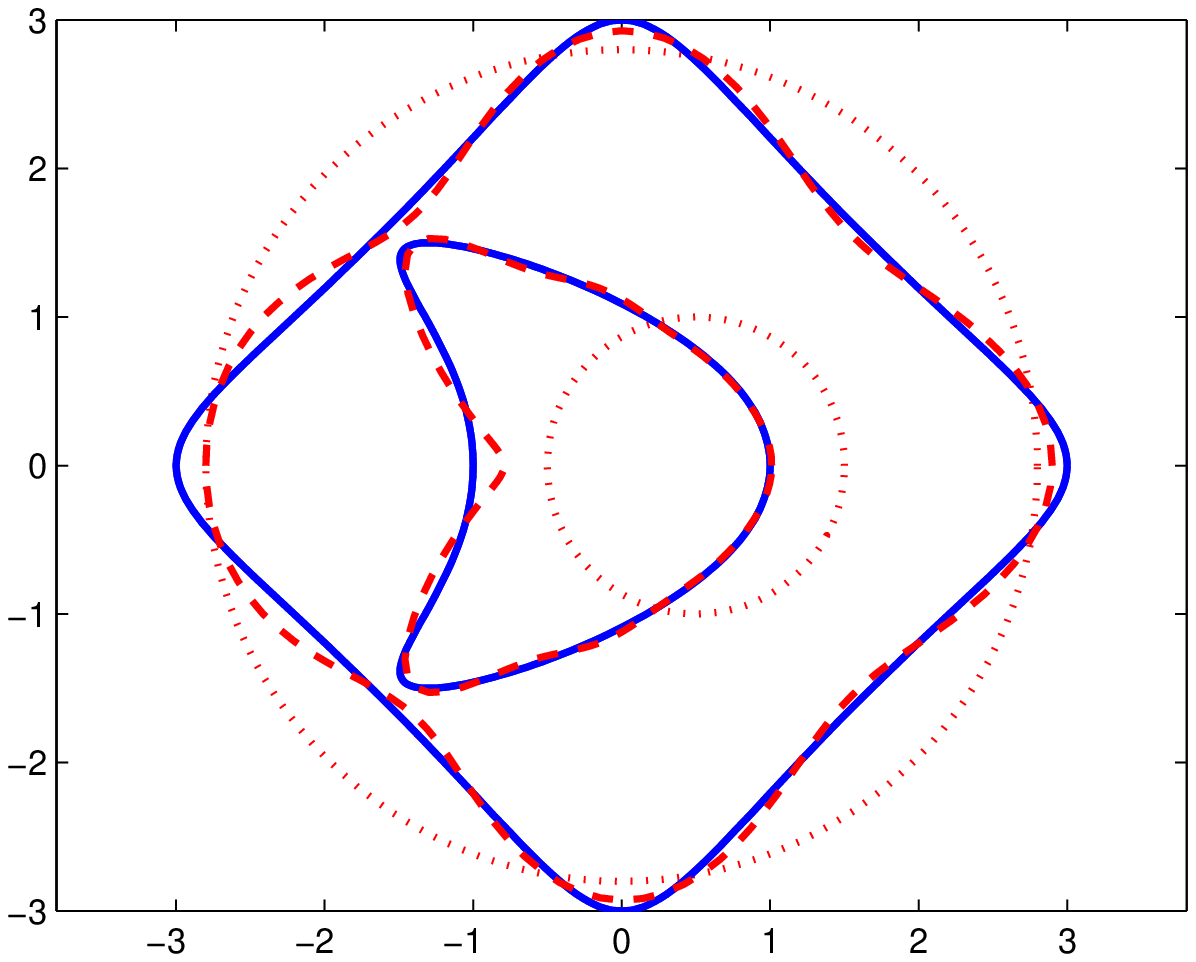}}
\caption{An kite-shaped, sound-hard obstacle embedded in a rounded square-shaped penetrable obstacle.
The boundary condition on $S_1$ is assumed unknown, and the exact data with a fixed wave number $k_0=1$
are used. (a), (b), (c) and (d) show the reconstruction results using three, four, eight and sixteen 
incident plane waves, respectively.}\label{fig3:subfig} 
\end{figure}

\begin{figure}[htbp]
\centering
\subfigure[\textbf{3\% noise}, $k_0=1$, $\la^{(Ir)}_1=-1.085\times10^{-3}$, $Err^{(Ir)}=0.04189$, 
$Ir=14$]{\label{fig3b:subfig:a}
\includegraphics[width=3in]{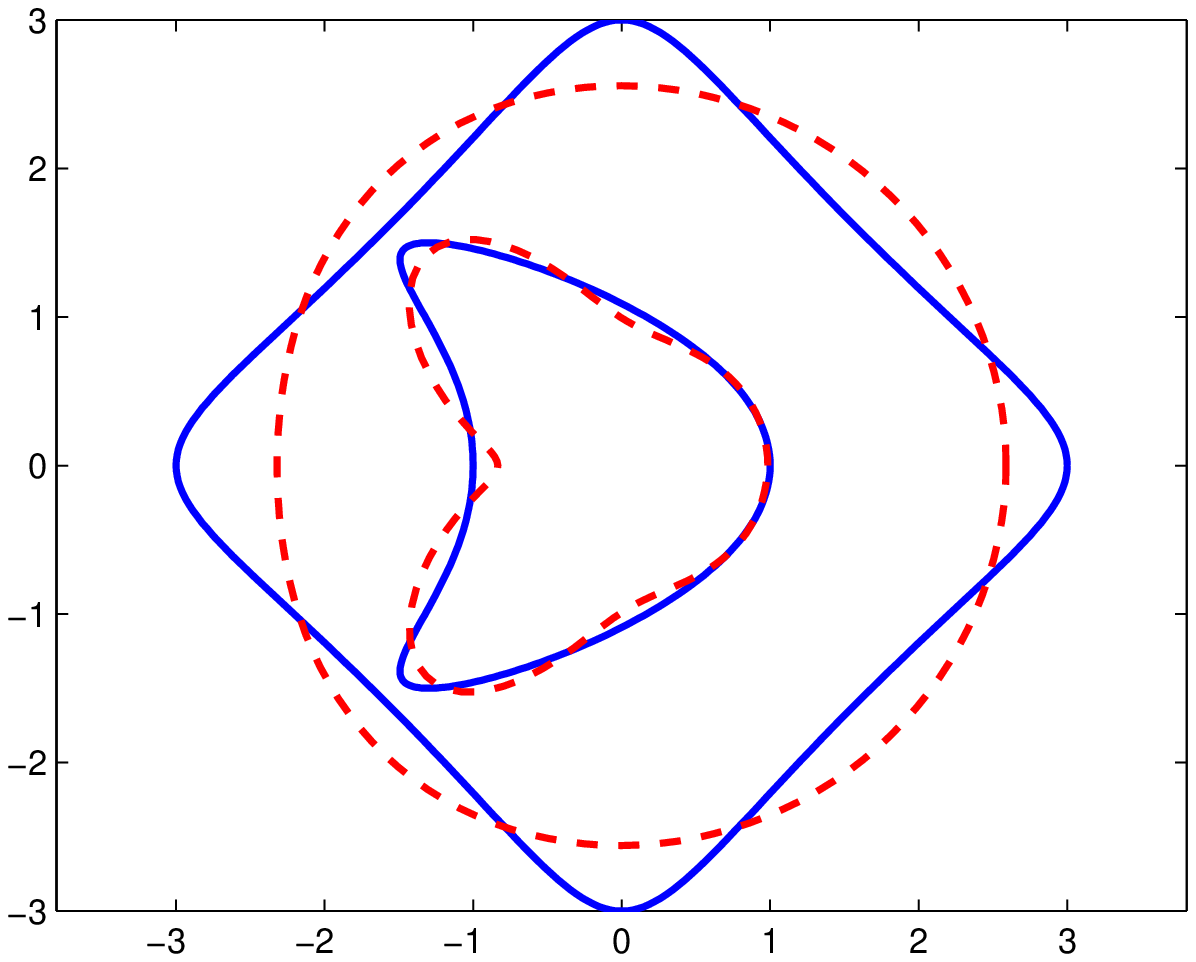}}
\hspace{0in}
\subfigure[\textbf{3\% noise}, $k_0=3$, $\la^{(Ir)}_1=-3.155\times10^{-3}$, $Err^{(Ir)}=0.04154$, 
$Ir=9$]{\label{fig3b:subfig:b}
\includegraphics[width=3in]{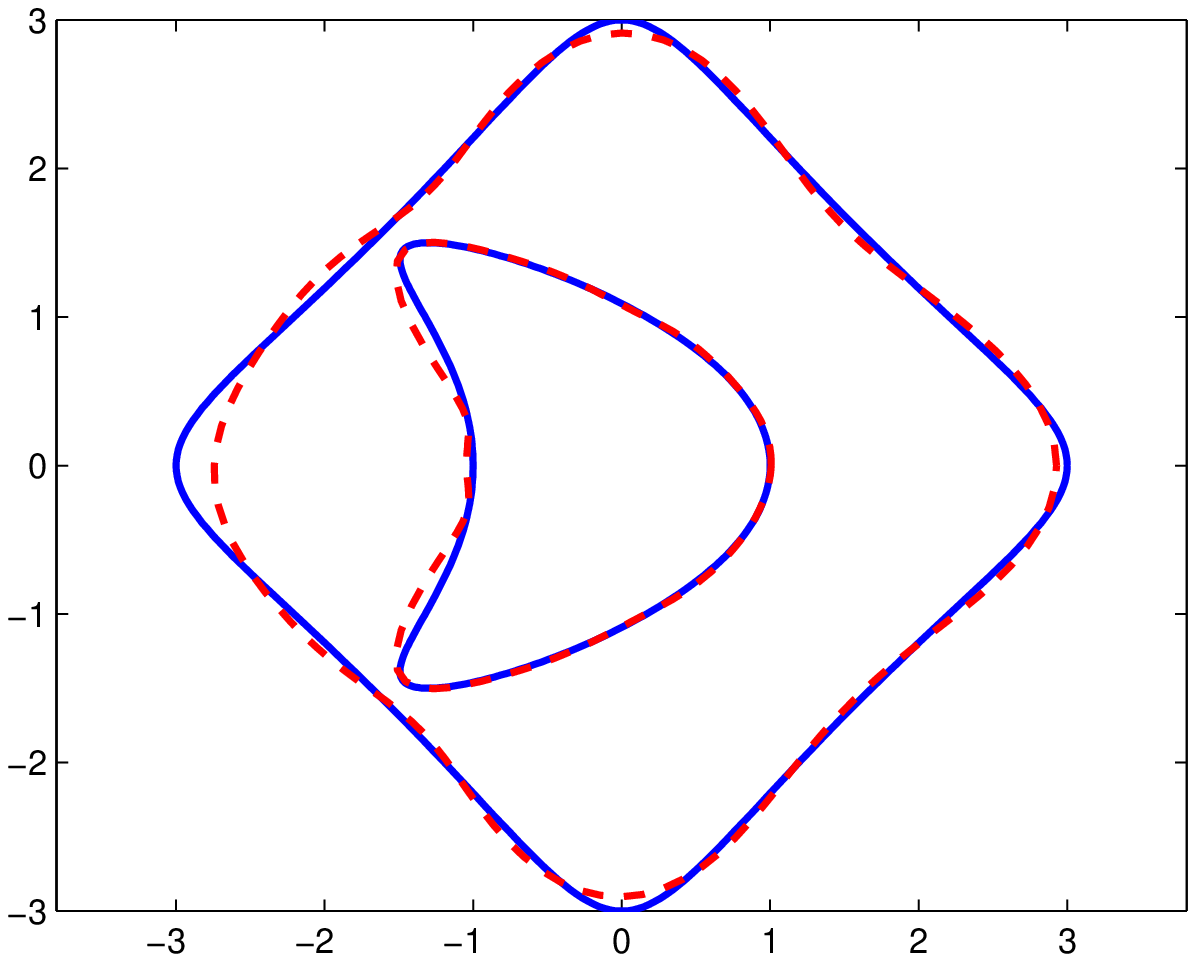}}
\subfigure[\textbf{3\% noise}, $k_0=5$, $\la^{(Ir)}_1=4.916\times10^{-3}$, $Err^{(Ir)}=0.04464$, 
$Ir=6$]{\label{fig3b:subfig:c} 
\includegraphics[width=3in]{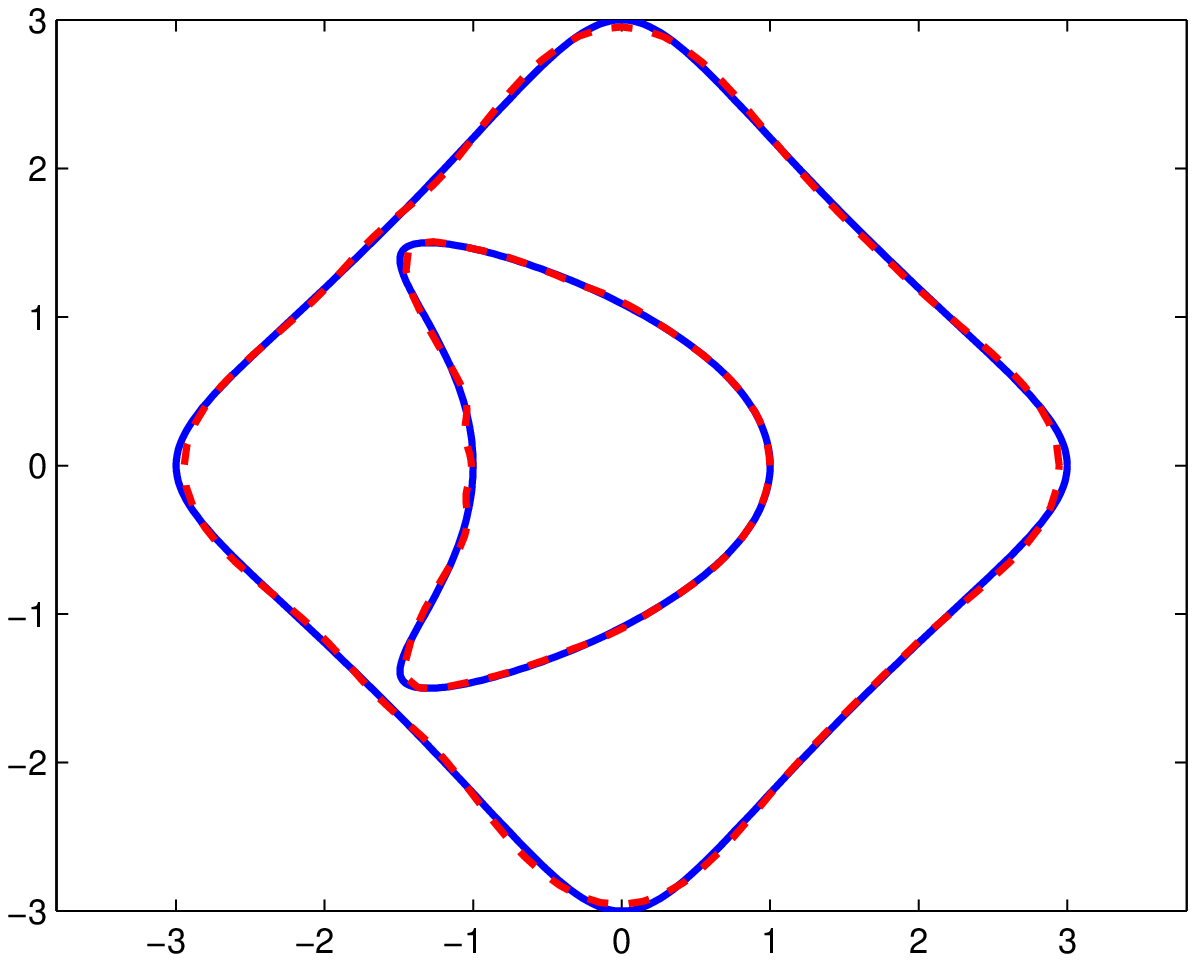}}
\hspace{0in}
\subfigure[\textbf{3\% noise}, $k_0=7$, $\la^{(Ir)}_1=-6.377\times10^{-3}$, $Err^{(Ir)}=0.04146$, 
$Ir=2$]{\label{fig3b:subfig:d} 
\includegraphics[width=3in]{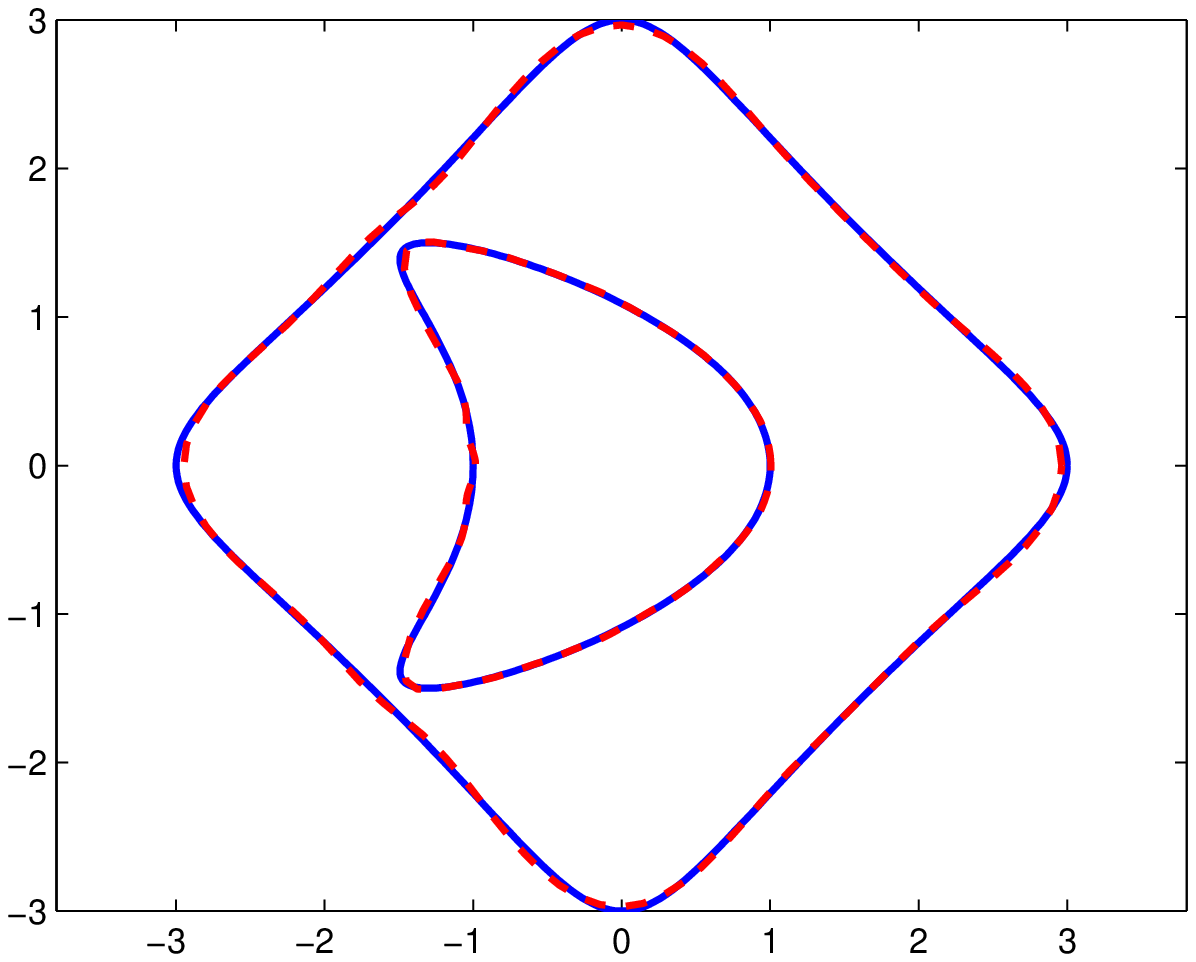}}
\caption{An kite-shaped, sound-hard obstacle embedded in a rounded square-shaped penetrable obstacle.
The boundary condition on $S_1$ is assumed unknown, and noisy multi-frequency data with
three incident plane waves are used. (a), (b), (c) and (d) show the reconstructions at $k_0=1,3,5,7,$
respectively.}\label{fig3b:subfig} 
\end{figure}

\textbf{Example 3.} In this example, we consider a sound-hard, kite-shaped obstacle
embedded in a rounded square-shaped penetrable obstacle.
We assume that the far-field data are given for this case which corresponds to a Neumann boundary
condition on the boundary $S_1$ of the embedded obstacle.
However, the boundary condition on $S_1$ is assumed unknown, and we will use our
method to reconstruct both obstacles and determine the type of boundary conditions on $S_1$.
For initial guesses of $S_0,S_1$, we choose the radii $r^{(0)}_0=2.8,r^{(0)}_1=1$ and 
the two centers to be $(0,0)$ and $(0.5,0)$, respectively. By trial, we choose $\la^{(0)}_1=1$. 
Since $S_0$ and $S_1$ are very close, it is more difficult to reconstruct both obstacles. 
In order to get a stable reconstruction, we use at least three incident plane waves. 
Moreover, the center of $\g^{(m)}_1$ is updated only at the first five iterations. 
Figure \ref{fig3:subfig} shows the reconstruction results using the exact far-field data 
with a fixed wave number $k_0=1$. Figures \ref{fig3:subfig:a}, \ref{fig3:subfig:b},
\ref{fig3:subfig:c} and \ref{fig3:subfig:d} give the reconstruction results with three, four, 
eight and sixteen incident plane waves, respectively. 
It is observed that using more incident plane waves does not yield a better reconstruction
even for the exact data.
In Figure \ref{fig3b:subfig} we present the reconstruction results using the $3\%$ noisy 
far-field data with three incident directions and multiple frequencies $k_0=1,3,5,7$. 
Compared with Figure \ref{fig3:subfig}, using multiple frequency data gives a much better 
reconstruction even for the noisy data. 
Since $|\la^{(Ir)}_1|<<1$ for each case, we conclude that $S_1$ is a sound-hard boundary.

\begin{figure}[htbp]
\centering
\subfigure[\textbf{3\% noise}, $k_0=1$, $Err^{(Ir)}=0.04124$, $Ir=13$]
{\label{fig4:subfig:a} 
\includegraphics[width=3in]{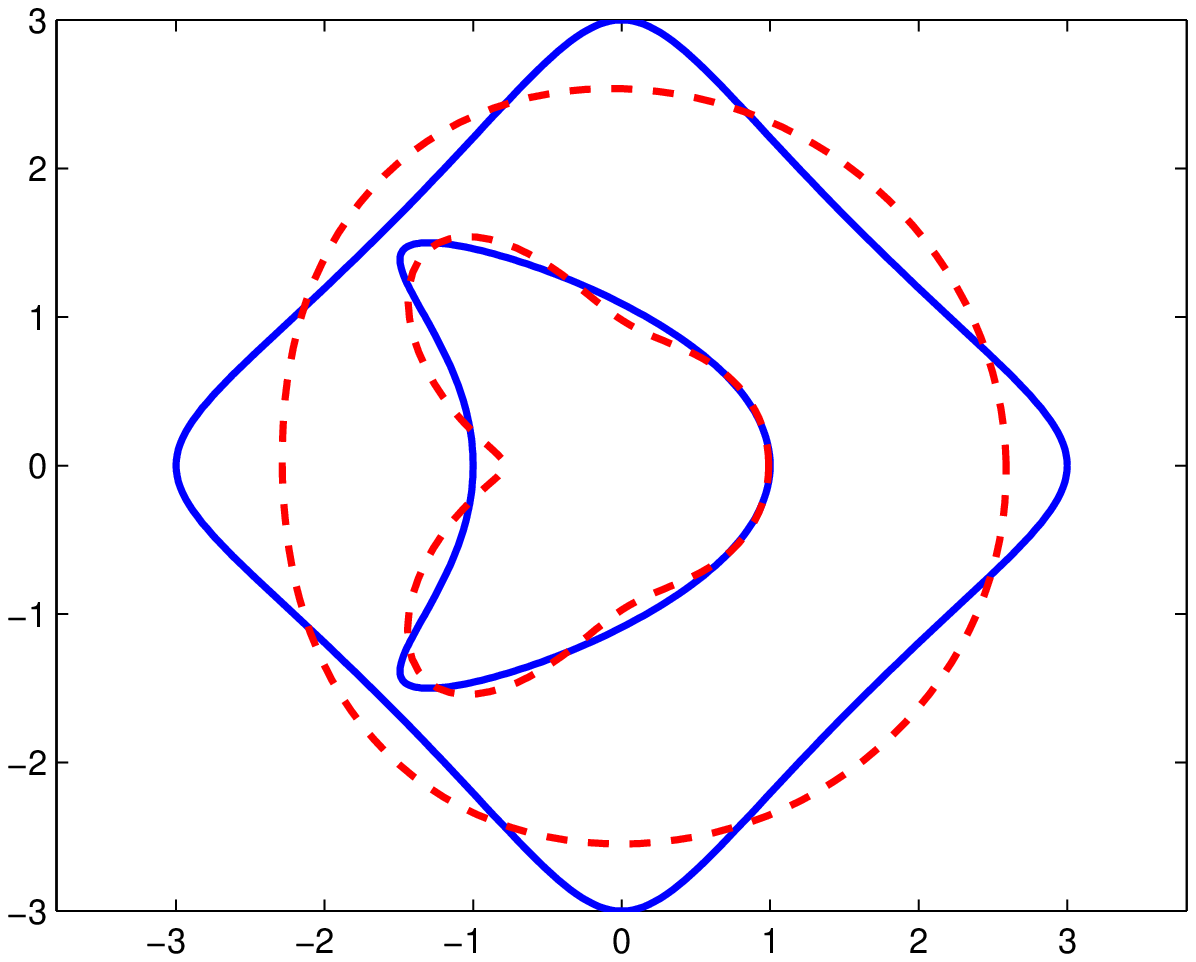}}
\hspace{0in}
\subfigure[\textbf{3\% noise}, $k_0=3$, $Err^{(Ir)}=0.03814$, $Ir=9$]
{\label{fig4:subfig:b} 
\includegraphics[width=3in]{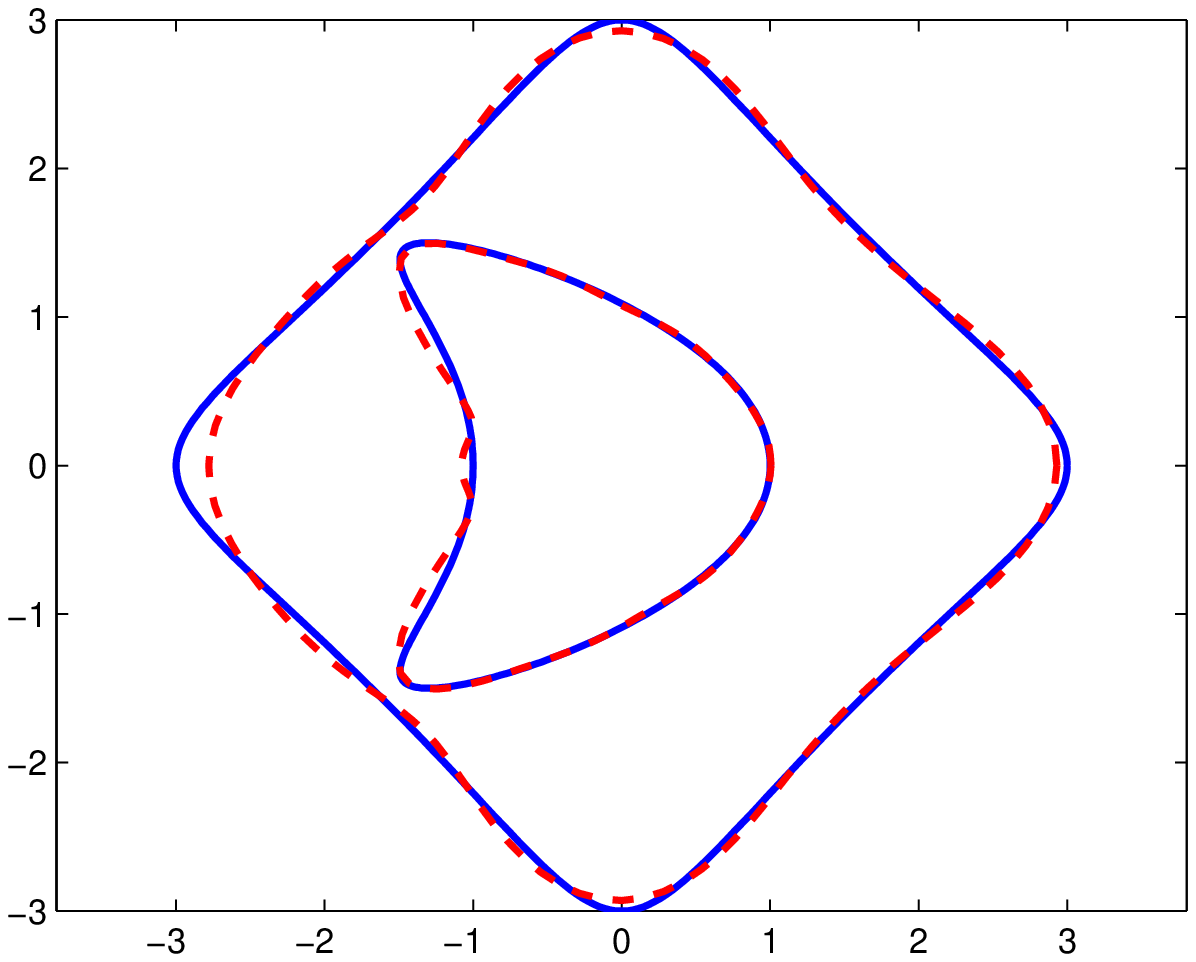}}
\subfigure[\textbf{3\% noise}, $k_0=5$, $Err^{(Ir)}=0.03943$, $Ir=6$]
{\label{fig4:subfig:c} 
\includegraphics[width=3in]{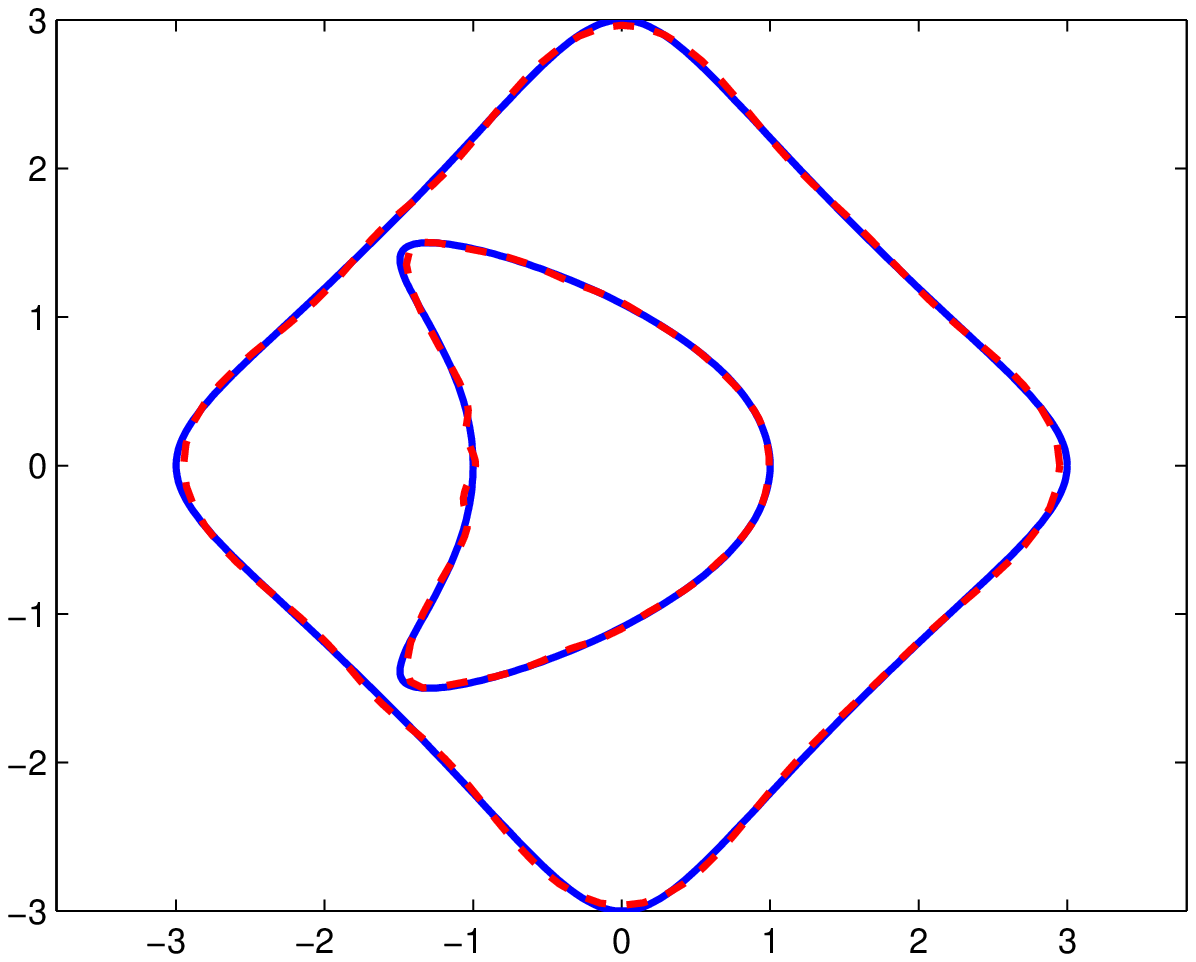}}
\hspace{0in}
\subfigure[\textbf{3\% noise}, $k_0=7$, $Err^{(Ir)}=0.04089$, $Ir=2$]
{\label{fig4:subfig:d} 
\includegraphics[width=3in]{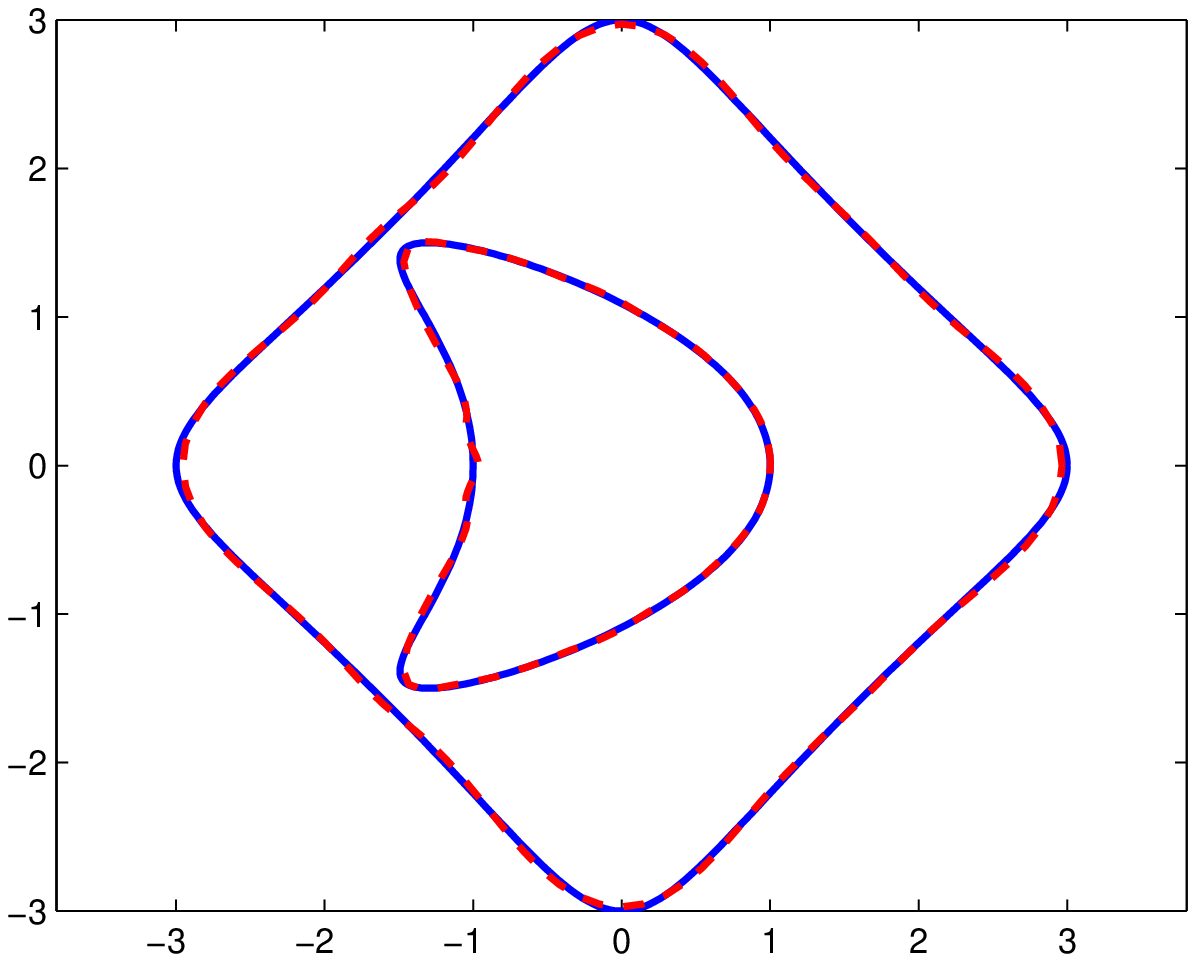}}
\caption{An kite-shaped, sound-hard obstacle embedded in a rounded square-shaped penetrable obstacle.
The boundary condition on $S_1$ is assumed known, and noisy multi-frequency data with
three incident plane waves are used. (a), (b), (c) and (d) show the reconstructions at $k_0=1,3,5,7,$ 
respectively.}\label{fig4:subfig} 
\end{figure}

\textbf{Example 4.} In this example, we consider the same problem as in Example 3.
Here, the far-field data are given again for the case of an embedded, sound-hard obstacle,
and $S_1$ is also assumed to be a sound-hard boundary.
We will use the Newton iteration method for the inverse problem corresponding to the original
scattering problem (\ref{eq16})-(\ref{rc}) with $\mathscr{B}(u):=\pa u/\pa\nu$ to reconstruct 
the two obstacles.
The initial guesses of $\g_0,\g_1,\la_0$ are the same as in Example 3, and the center of
$\g^{(m)}_1$ is updated only at the first five iterations.
In Figure \ref{fig4:subfig}, we present the reconstruction results by using
the $3\%$ noisy far-field data with three incident direction and multiple frequencies $k_0=2,4,6,8$.
Figures \ref{fig3b:subfig} and \ref{fig4:subfig} indicate that the reconstruction obtained by both 
methods is almost the same.

\appendix
\renewcommand{\theequation}{\Alph{section}.\arabic{equation}}
\section*{Appendix: Frechet derivatives}\label{sec-a}
\setcounter{section}{1}
\setcounter{equation}{0}

In this appendix we study the Frechet derivative of the far field operators $F$ and $\tilde{F}$ 
when $\la_1>0$ and $\tau_1>0,$ respectively.

Let $(h_0,h_1,\Delta\la_1)\in C^1(S_0)\times C^1(S_1)\times\R$ be small perturbations
of $(S_0,S_1,\la_1)$ (that is, $||h_i||_{C^1(S_i)},\;\Delta\la_1<<1,\;i=0,1$).
Define $(S_i)_{h_i}\coloneqq\{y\in\R^2\,|\,y=x+h_i(x),\,x\in S_i\},\;i=0,1$.
The far field operator $F$ is called Frechet differentiable at $(S_0,S_1,\la_1)$ if there exists 
a linear bounded operator
$F'|_{(S_0,S_1,\la_1)}:C^1(S_0)\times C^1(S_1)\times\R\rightarrow L^2(S^1)$ such that
\ben
&&\|F((S_0)_{h_0},(S_1)_{h_1},\la_1+\Delta\la_1)-F(S_0,S_1,\la_1)
-F'|_{(S_0,S_1,\la_1)}(h_0,h_1,\Delta\la_1)\|_{L^2(S^1)}\\
&&\qquad\qquad=o\left(||h_0||_{C^1(S_0)}+||h_1||_{C^1(S_1)}+|\Delta\la_1|\right)
\enn
$F'|_{(S_0,S_1,\la_1)}$ is called the Frechet derivative of $F$ at $(S_0,S_1,\la_1)$. 
Following the idea in \cite{Hettlich1995,Kirsch93}, we have the following theorem to characterize  
the Frechet derivation of $F$ when $\la_1>0$.

\begin{theorem}\label{th3}
Let $S_0,S_1$ be $C^2$ and let $u\in H^1_{loc}(\R^2)$ be 
the solution of the problem (\ref{eq20})-(\ref{eq24}) and (\ref{rc}) with $\la_1>0.$
Then $F$ is Frechet differentiable at $(S_0,S_1,\la_1)$ and the Frechet derivative 
$F'|_{(S_0,S_1,\la_1)}$ is given by $F'|_{(S_0,S_1,\la_1)}(h_0,h_1,\Delta\la_1)=u'_\infty$, 
where $u'_\infty$ is the far-field pattern of 
$u'\in H^1_{loc}(\Om_0)\bigcap H^1(\Om_1)\bigcap H^1(\Om_2)$ which solves the problem 
(\ref{eq39})-(\ref{eq11}) and (\ref{rc}) with 
\ben
f_1&=&-(h_0)_\nu\left(\frac{\pa u_+}{\pa\nu}-\frac{\pa u_-}{\pa\nu}\right)\Big|_{S_0},\\
f_2&=&(k^2_0-\la_0 k^2_1)(h_0)_\nu u\big|_{S_0}+Div_{S_0}\left[(h_0)_\nu\big((\na u_+)_t
-\la_0(\na u_-)_t\big)\right],\\
f_3&=&-(h_1)_\nu\left(\frac{\pa u_+}{\pa\nu}-\frac{\pa u_-}{\pa\nu}\right)\Big|_{S_1},\\
f_4&=&(k^2_1-\la_1 k^2_2)(h_1)_\nu u\big|_{S_1}+Div_{S_1}\left[(h_1)_\nu\big((\na u_+)_t
-\la_1(\na u_-)_t\big)\right]+\frac{\Delta\la_1}{\la_1}\frac{\pa u_+}{\pa\nu}\Big|_{S_1},
\enn
where $h_\nu$ and $h_t$ denote the normal and tangential components of the vector field $h$.
\end{theorem}

\begin{remark}\label{re-a} {\rm 
From \cite{Peter2000} it follows that $u\in H^2_{loc}(\Om_0)\bigcap H^2(\Om_1)\bigcap H^2(\Om_2),$
which guarantees the existence of $u'$.
}
\end{remark}

\begin{proof}
Our proof is based on the variational formulation. To this end, let $B_R$ be a fixed ball of
radius $R$ and centered at the original with $R$ large enough such that 
$\overline{\Om_1\bigcup\Om_2}\subset B_R$ and let $L$ be the Dirichlet to Neumann (D-N) map 
defined on $\pa B_R$ \cite{ColtonKress98}. It is easy to see that the solution $u$ satisfies that
for all $v\in H^1(B_R)$,
\be\label{eq54} 
\int_{B_R}(\al\nabla u\cdot\nabla\ol{v}-k^2 u\ol{v})dx-<Lu,v>
=\int_{\pa B_R}\left(\frac{\pa u^i}{\pa\nu}-L u^i\right)\ol{v}dx
\en
with
\ben
\alpha=\left\{\begin{aligned}1\quad&\mbox{in}\;\;\Om_0,\\
                    \la_0\quad&\mbox{in}\;\;\Om_1,\\
                    \la_0\la_1\quad&\mbox{in}\;\;\Om_2,
               \end{aligned}\right.\qquad
k^2=\left\{\begin{aligned}k^2_0\quad&\mbox{in}\;\;\Om_0,\\
              \la_0k^2_1\quad&\mbox{in}\;\;\Om_1,\\
              \la_0\la_1 k^2_2\quad&\mbox{in}\;\;\Om_2,
           \end{aligned}\right.
\enn
where $<\cdot,\cdot>$ denotes the dual between $H^{-1/2}(\pa B_R)$ and $H^{1/2}(\pa B_R)$.
For the property of $L$ and the proof of the equivalence between the variational equation and 
the origin problem, we refer to \cite{LiuZhangHu2010} where a similar problem is considered.

Denote by $S(u,v)$ the left-hand side of (\ref{eq54}). Then, by the Fredholm alternative and 
the uniqueness of the origin problem, it can be proved that there exists a linear bounded operator
$B:H^1(B_R)\rightarrow H^1(B_R)$ with a bounded inverse such that 
\be\label{eq55} 
S(u,v)=(B(u),v),
\en 
where $(\cdot,\cdot)$ denotes the inner product in $H^{1}(B_R)$. 
Thus, for any $f\in H^1(B_R)$ there exists a unique function $u_f\in H^1(B_R)$ such that
$S(u_f,v)=(f,v)$ for all $v\in H^1(B_R)$.

Let $K_0,K_1$ be two domains with $S_0\subset K_0\subset\subset B_R\ba\ol{\Om}_2$, 
$S_1\subset K_1\subset\subset B_R\ba\ol{\Om}_0$ and $\ol{K}_0\cap\ol{K}_1=\emptyset$. 
For $(h_0,h_1)\in C^1(S_0)\times C^1(S_1)$ there exists an extension of $(h_0,h_1)$ 
which is denoted again by $(h_0,h_1)$ such that $h_0,h_1\in C^1(\mathbb{R}^2)$,
$\textrm{supp}(h_0)\subset K_0$, $\textrm{supp}(h_1)\subset K_1$ and
$||h_i||_{C^1(\R^2)}\leq C||h_i||_{C^1(S_i)},\;i=0,1$. 
Let $h(x)=h_0(x)+h_1(x)$ for $x\in\mathbb{R}^2$ and let $\varphi(x)=x+h(x)$.
If $||h_i||_{C^1(S_i)},\;i=0,1,$ are small enough, then $\varphi:B_R\rightarrow B_R$ is a 
diffeomorphism (see \cite{Hohage99}) and $\varphi$ maps $S_0$ to $(S_0)_{h_0}$, 
$S_1$ to $(S_1)_{h_1}$, respectively. Denote by $\psi$ the inverse of $\varphi$ and 
by $J_\varphi,\;J_\psi$ the Jacobians of $\varphi,\psi$, respectively.

Let $u_{(h,\Delta\la_1)}$ be the solution of the problem (\ref{eq20})-(\ref{eq24}) and (\ref{rc}) 
with ($S_0,S_1,\lambda_1$) replaced by $((S_0)_{h_0},(S_1)_{h_1},\la_1+\Delta\la_1)$.
Then $u_{(h,\Delta\la_1)}$ satisfies that for all $v\in{H}^1(B_R)$,
\be\label{eq56}
\int_{B_R}(\al_{(h,\Delta\la_1)}\nabla u_{(h,\Delta\la_1)}\cdot\nabla\ol{v}
-k^2_{(h,\Delta\la_1)}u_{(h,\Delta\la_1)}\ol{v})dx-<Lu_{(h,\Delta\la_1)},v>
=\int_{\pa B_R}\left(\frac{\pa u^i}{\pa\nu}-L u^i\right)\ol{v}dx
\en
with
\ben
\alpha_{(h,\Delta\la_1)}=\left\{\begin{aligned}1\quad&\mbox{in}\;\;(\Om_0)_h,\\
             \la_0\quad&\mbox{in}\;\;(\Om_1)_h,\\
             \la_0(\la_1+\triangle\la_1)\quad&\mbox{in}\;\;(\Om_2)_h,
             \end{aligned}\right.\qquad
k^2_{(h,\Delta\la_1)}=\left\{\begin{aligned}k^2_0\quad&\mbox{in}\;\;(\Om_0)_h,\\
                                            \la_0k^2_1\quad&\mbox{in}\;\;(\Om_1)_h,\\
                                            \la_0(\la_1+\Delta\la_1)k^2_2\quad&\mbox{in}\;\;(\Om_2)_h,
                             \end{aligned}\right.
\enn
where $(\Om_i)_h=\{x+h(x)\,|\,x\in\Om_i\},\;i=0,1,2$. 
If $\Delta\la_1$ is small enough, then the unique solution of the problem (\ref{eq56}) 
exists and can be extended outside of $B_R$ such that $u_{(h,\Delta\la_1)}\in H^1_{loc}(\R^2)$ 
with $u_{(h,\Delta\la_1)}$ satisfying (\ref{eq20}) and $u_{(h,\Delta\la_1)}-u^i$ satisfying (\ref{rc}).

Write $\tilde{u}_{(h,\Delta\la_1)}=u_{(h,\Delta\la_1)}\circ\varphi$. 
Then, by the change of variables, we obtain that for all $v\in H^1(B_R),$
\ben
&&\int_{B_R}(\al\nabla\tilde{u}^T_{(h,\Delta\la_1)}J_\psi J_\psi^T\nabla\ol{v}
-k^2\tilde{u}_{(h,\Delta\la_1)}\ol{v})\det(J_\varphi)dx\\
&&\quad+\la_0\Delta\la_1\int_{\Om_2}(\nabla\tilde{u}^T_{(h,\Delta\la_1)}J_\psi J_\psi^T\nabla\ol{v}
  -k^2_2\tilde{u}_{(h,\Delta\la_1)}\ol{v})\det(J_\varphi)dx
  -<L\tilde{u}_{(h,\Delta\la_1)},v>
=\int_{\pa B_R}\left(\frac{\pa u^i}{\pa\nu}-Lu^i\right)\ol{v}dx
\enn
Denote by $S_{(h,\Delta\la_1)}(\tilde{u}_{(h,\Delta\la_1)},v)$ the left-hand side of the above 
equation. Then, for all $v\in H^1(B_R)$ we have
\ben
S(u-\tilde{u}_{(h,\Delta\la_1)},v)&=&S_{(h,\Delta\la_1)}(\tilde{u}_{(h,\Delta\la_1)},v)
  -S(\tilde{u}_{(h,\Delta\la_1)},v)\\
&=&\int_{B_R}\left[\al\nabla\tilde{u}^T_{(h,\Delta\la_1)}
   \left(J_\psi J_\psi^T\det(J_\varphi)-I\right)\nabla\ol{v}
   -k^2\left(\det(J_\varphi)-1\right)\tilde{u}_{(h,\Delta\la_1)}\ol{v}\right]dx\\
&&+\la_0\Delta\la_1\int_{\Om_2}\left(\nabla\tilde{u}^T_{(h,\Delta\la_1)}J_\psi J_\psi^T\nabla\ol{v}
  -k^2_2\tilde{u}_{(h,\Delta\la_1)}\ol{v}\right)\det(J_\varphi)dx
\enn

Let $w\in H^1(B_R)$ be the solution of the problem: find $w\in H^1(B_R)$ such that 
\be\label{eq1}
S(w,v)&=&\int_{B_R}\left[\al\na u^T\left(J_h+J_h^T-\na\cdot(hI)\right)\nabla\ol{v}
         +k^2(\nabla\cdot h)u\ol{v}\right]dx\no\\
&&-\la_0\Delta\la_1\int_{\Om_2}\left(\nabla u^T\nabla\ol{v}-k^2_2 u\ol{v}\right)dx,
   \qquad\forall v\in H^1(B_R)
\en
(by the property of $S(\cdot,\cdot)$, it is easy to prove that such a $w$ exists).
Then, and since $|J_\psi J_\psi^T\det(J_\varphi)-I+(J_h+J_h^T-(\nabla\cdot h)I)|
=O\left(||h||^2_{C^1(\ol{B_R})}\right)$, 
$\det(J_\varphi)-1=\nabla\cdot h+O\left(||h||^2_{C^1(\ol{B_R})}\right)$ and
$\tilde{u}_{(h,\Delta\la_1)}$ tends to $u$ in $H^1(B_R)$ (which can be proved similarly 
as in \cite{Hettlich1995}), it is easy to prove that 
\be\label{eq56+}
\sup\limits_{v\in H^1(B_R)}\frac{S(\tilde{u}_{(h,\Delta\la_1)}-u-w,v)}{||v||_{H^1(B_R)}}
=o\left(||h||_{C^1(\ol{B_R})}+|\Delta\la_1|\right).
\en
Then, by (\ref{eq55}) we have 
$[{1}/({||h||_{C^1(B_R)}+|\Delta\la_1|})](\tilde{u}_{(h,\Delta\la_1)}-u-w)\rightarrow 0$ 
in $H^1(B_R)$ as $||h||_{C^1(B_R)}+|\Delta\la_1|$ tends to zero. 
This, together with the trace theorem, implies that 
$[{1}/({||h||_{C^1(B_R)}+|\Delta\la_1|})](\tilde{u}_{(h,\Delta\la_1)}-u-w)\rightarrow 0$
in $H^{1/2}(\pa B_R)$ as $||h||_{C^1(B_R)}+|\Delta\la_1|$ tends to zero.

Now let $u'=w-h\cdot\nabla u$. Then $u'\in H^1(B_R\cap\Om_0)\bigcap H^1(\Om_1)\bigcap H^1(\Om_2)$
since $u\in H^2_{loc}(\Om_0)\bigcap H^2(\Om_1)\bigcap H^2(\Om_2)$.
Note that 
$$
\nabla u^T[J_h+J_h^T-(\nabla\cdot h)I]\nabla\ol{v}=\nabla\cdot[(h\cdot\nabla\ol{v})\nabla u
-(\nabla u\cdot\nabla\ol{v})h]+\nabla(h\cdot\nabla u)\cdot\nabla\ol{v}-(h\cdot\nabla\ol{v})\Delta u.
$$
Then for $v\in H^2(B_R)$ the right-hand side of (\ref{eq1}) becomes
\ben
&&\int_{B_R}\left[\al\nabla(h\cdot\nabla u)\cdot\nabla\ol{v}
  +\al\nabla\cdot\big((h\cdot\nabla\ol{v})\nabla u-(\nabla u\cdot\nabla\ol{v})h\big)
  -\al(h\cdot\nabla\ol{v})\Delta u+k^2(\nabla\cdot h)u\ol{v}\right]dx\\
&&\qquad\quad-\la_0\Delta\la_1\int_{\Om_2}(\nabla u^T\nabla\ol{v}-k^2_2u\ol{v})dx\\
&&\quad=\int_{B_R}\left[\al\nabla(h\cdot\nabla u)\cdot\nabla\ol{v}-k^2(h\cdot\nabla u)\ol{v}
   +\al\nabla\cdot\big((h\cdot\nabla\ol{v})\nabla u-(\nabla u\cdot\nabla\ol{v})h\big)
   +k^2\nabla\cdot(hu\ol{v})\right]dx\\
&&\qquad\quad-\la_0\Delta\la_1\int_{\Om_2}\nabla\cdot(\nabla u\ol{v})dx.
\enn
On the other hand, the left-hand side of (\ref{eq1}) is 
$S(w,v)=\int_{B_R}(\al\nabla w\cdot\na\ol{v}-k^2w\ol{v})dx-<Lw,v>$. 
Thus, we have that for $v\in H^2(B_R)$,
\ben
S(u',v)=\int_{B_R}\left[\al\nabla\cdot\left((h\cdot\nabla\ol{v})\nabla u
 -(\nabla u\cdot\nabla\ol{v})h\right)+k^2\nabla\cdot(hu\ol{v})\right]dx
 -\la_0\Delta\la_1\int_{\Om_2}\nabla\cdot(\nabla u\ol{v})dx
\enn 
where we have used the fact that $h(x)=0$ in the neighborhood of $\pa B_R.$
By the divergence theorem, the right-hand side of the above equation becomes
\be\no
&&-\int_{S_0}\left[\left((h\cdot\nabla\ol{v})(\nabla u)_+
  -((\nabla u)_+\cdot\nabla\ol{v})h\right)\cdot\nu+k^2_0h_\nu u\ol{v}\right]ds\\ \no
&&+\la_0\int_{S_0}\left[\left((h\cdot\nabla\ol{v})(\nabla u)_-
  -((\nabla u)_-\cdot\nabla\ol{v})h\right)\cdot\nu+k^2_1h_\nu u\ol{v}\right]ds\\ \no
&&-\la_0\int_{S_1}\left[\left((h\cdot\nabla\ol{v})(\nabla u)_+
  -((\nabla u)_+\cdot\nabla\ol{v})h\right)\cdot\nu+k^2_1h_\nu u\ol{v}\right]ds\\ \no
&&+\la_0\la_1\int_{S_1}\left[\left((h\cdot\nabla\ol{v})(\nabla u)_-
  -((\nabla u)_-\cdot\nabla\ol{v})h\right)\cdot\nu+k^2_2h_\nu u\ol{v}\right]ds
  -\la_0\Delta\la_1\int_{S_1}\frac{\pa u_-}{\pa\nu}\ol{v}ds\\ \no
&=&-\int_{S_0}\left[(h\cdot\nabla\ol{v})\left(\frac{\pa u_+}{\pa\nu}
   -\la_0\frac{\pa u_-}{\pa\nu}\right)-h_\nu((\nabla u)_+
   -\la_0(\nabla u)_-)\cdot\nabla\ol{v}+(k^2_0-\la_0k^2_1)h_\nu u\ol{v}\right]ds\\ \no
&&-\la_0\int_{S_1}\left[(h\cdot\nabla\ol{v})\left(\frac{\pa u_+}{\pa\nu}
  -\la_1\frac{\pa u_-}{\pa\nu}\right)-h_\nu((\nabla u)_+
  -\la_1(\nabla u)_-)\cdot\nabla\ol{v}+(k^2_1-\la_1 k^2_2)h_\nu u\ol{v}\right]ds\\ \no
&&-\la_0\Delta\la_1\int_{S_1}\frac{\pa u_-}{\pa\nu}\ol{v}ds\\ \no
&=&\int_{S_0}\left[h_\nu((\nabla u_+)_t-\la_0(\nabla u_-)_t)\cdot(\nabla\ol{v})_t
   -(k^2_0-\la_0 k^2_1)h_\nu u\ol{v}\right]ds\\ \label{eq57}
&&+\la_0\int_{S_1}\left[h_\nu((\nabla u_+)_t-\la_1(\nabla u_-)_t)\cdot(\nabla\ol{v})_t
  -(k^2_1-\la_1 k^2_2)h_\nu u\ol{v}\right]ds
  -\la_0\Delta\la_1\int_{S_1}\frac{\pa u_-}{\pa\nu}\ol{v}ds,
\en
where we have used the transmission conditions (\ref{eq23}) and (\ref{eq24}) on $S_i,\;i=0,1.$
Then we conclude that $u'$ satisfies the equation (\ref{eq39}) in $B_R\bigcap\Om_0$ and the equations
(\ref{eq39}) and (\ref{eq10}) and that $Lu'={\pa u'_-}/{\pa\nu}\big|_{\pa B_R}$. 
We extend $u'$ outside of $B_R$ to be the solution of the Helmholtz equation $\Delta U+k_0^2U=0$  
satisfying the Dirichlet boundary condition $U=u'_-$ on $\pa B_R$ and 
the radiation condition (\ref{rc}). Denote the extension of $u'$ by $u'$ again.
Then $u'$ satisfies (\ref{eq39}) in $\Om_0$, and we have
\ben
S(u',v)=-\int_{S_0}\left(\frac{\pa u'_+}{\pa\nu}-\la_0\frac{\pa u'_-}{\pa\nu}\right)\ol{v}ds
        -\la_0\int_{S_1}\left(\frac{\pa u'_+}{\pa\nu}-\la_1\frac{\pa u'_-}{\pa\nu}\right)\ol{v}ds.
\enn
From this and (\ref{eq57}) it follows that 
\ben
\frac{\pa u'_+}{\pa\nu}-\la_0\frac{\pa u'_-}{\pa\nu}&=&(k^2_0-\la_0 k^2_1)(h_0)_\nu u\big|_{S_0}
  +Div_{S_0}\Big[(h_0)_\nu\big((\nabla u_+)_t-\la_0(\nabla u_-)_t\big)\Big]\quad\mbox{on}\;\;S_0,\\
\frac{\pa u'_+}{\pa\nu}-\la_1\frac{\pa u'_-}{\pa\nu}&=&(k^2_1-\la_1 k^2_2)(h_1)_\nu u\big|_{S_1}
  +Div_{S_1}\Big[(h_1)_\nu\big((\nabla u_+)_t-\la_1(\nabla u_-)_t\big)\Big]
  +\frac{\Delta\la_1}{\la_1}\frac{\pa u_+}{\pa\nu}\Big|_{S_1}\quad\mbox{on}\;\;S_1.
\enn
The definition of $u'$ gives 
\ben
u'_+-u'_-&=&-(h_0)_\nu\left(\frac{\pa u_+}{\pa\nu}
            -\frac{\pa u_-}{\pa\nu}\right)\Big|_{S_0}\quad\mbox{on}\;\;S_0,\\
u'_+-u'_-&=&-(h_1)_\nu\left(\frac{\pa u_+}{\pa\nu}
            -\frac{\pa u_-}{\pa\nu}\right)\Big|_{S_1}\quad\mbox{on}\;\;S_1.\\
\enn

Since $supp(h)\subset B_R$, we have $\tilde{u}_{(h,\Delta\la_1)}-u-w={u}_{(h,\Delta\la_1)}-u-u'$
on $\pa B_R$. This, together with the discussion just following (\ref{eq56+}), implies that 
$({u}^s_{(h,\Delta\la_1)}-u^s-u')_\infty=o(||h_0||_{C^1(S_0)}+||h_1||_{C^1(S_1)}+|\Delta\la_1|)$
as $||h_0||_{C^1(S_0)}+||h_1||_{C^1(S_1)}+|\Delta\la_1|\rightarrow0$.
We thus conclude that $F$ is Frechet differentiable at $(S_0,S_1,\la_1)$ and 
the Frechet derivative $F'|_{(S_0,S_1,\la_1)}(h_0,h_1,\Delta\la_1)=u'_\infty$.
The proof is then completed.
\end{proof}

The following theorem can be shown similarly.

\begin{theorem}\label{th4}
Let $S_0,S_1$ be $C^2$ and let $u\in H^1_{loc}(\R^2)$ be the solution of the 
problem (\ref{eq20})-(\ref{eq21}), (\ref{eq23}), (\ref{eq22a})-(\ref{eq24a}) and (\ref{rc})
with $\tau_1>0$. Then $\tilde{F}$ is Frechet differentiable at $(S_0,S_1,\tau_1)$ and
the Frechet derivative $\tilde{F}'|_{(S_0,S_1,\tau_1)}$ is given by 
$\tilde{F}'|_{(S_0,S_1,\tau_1)}(h_0,h_1,\Delta\tau_1)=u'_\infty$,
where $u'_\infty$ is the far-field pattern of 
$u'\in H^1_{loc}(\Om_0)\bigcap H^1(\Om_1)\bigcap H^1(\Om_2)$ which solves the problem 
(\ref{eq39})-(\ref{eq11}) and (\ref{rc}) with
\ben
f_1&=&-(h_0)_\nu\left(\frac{\pa u_+}{\pa\nu}-\frac{\pa u_-}{\pa\nu}\right)\Big|_{S_0},\\
f_2&=&(k^2_0-\la_0 k^2_1)(h_0)_\nu u\big|_{S_0}
      +Div_{S_0}\Big[(h_0)_\nu\big((\nabla u_+)_t-\la_0(\nabla u_-)_t\big)\Big],\\
f_3&=&-(h_1)_\nu\left(\frac{\pa u_+}{\pa\nu}-\frac{\pa u_-}{\pa\nu}\right)\Big|_{S_1},\\
f_4&=&(k^2_1-\la_1 k^2_2)(h_1)_\nu u\big|_{S_1}
      +Div_{S_1}\Big[(h_1)_\nu\big((\nabla u_+)_t-\la_1(\nabla u_-)_t\big)\Big]
      -\frac{\triangle\tau_1}{\tau_1}\frac{\pa u_+}{\pa\nu}\Big|_{S_1}.
\enn
\end{theorem}

\begin{remark}\label{re9} {\rm
By \cite{Peter2000} again, $u\in H^2_{loc}(\Om_0)\bigcap H^2(\Om_1)\bigcap H^2(\Om_2)$, 
so $u'$ exists.
}
\end{remark}

\section*{Acknowledgements}
This work was supported by the NNSF of China under grants 11071244 and 11161130002.


\end{document}